\documentclass[preprint,12pt]{elsarticle}

\usepackage{amssymb}
\usepackage{amsmath}
\usepackage{amsthm}
\usepackage{lmodern}
\usepackage{enumerate}

\usepackage{tikz}
\usetikzlibrary{calc,arrows.meta,decorations.pathreplacing}

\newtheorem{thm}{Theorem}[section]
\newtheorem{prop}[thm]{Proposition}
\newtheorem{lem}[thm]{Lemma}

\theoremstyle{definition}
\newtheorem{df}[thm]{Definition}
\newtheorem{example}[thm]{Example}
\newtheorem{as}[thm]{Assumption}

\theoremstyle{remark}

\newcommand{\bfx}{{ x}}

\newcommand{\xbf}{{ x}}

\newcommand{\ybf}{{ y}}
\newcommand{\zbf}{{ z}}
\newcommand{\Xbf}{{ X}}

\newcommand{\Ybf}{{ Y}}
\newcommand{\Zbf}{{ Z}}
\newcommand{\Wbf}{{ W}}

\newcommand{\Mbf}{{ M}}
\newcommand{\Ubf}{{ U}}
\newcommand{\Nbf}{{ N}}
\newcommand{\Sbf}{{ S}}

\newcommand{\Tbf}{{ T}}
\newcommand{\Rbf}{{ R}}
\newcommand{\mbf}{{ m}}
\newcommand{\nbf}{{ n}}

\newcommand{\bfy}{{ y}}

\newcommand{\eps}{{\varepsilon}}

\newcommand{\1}{{\mathbf 1}}
\renewcommand{\P}{\mathbb{P}}

\numberwithin{equation}{section}

\newcommand{\Ebb}{\mathbb{E}}  
 
\newcommand{\E}{\mathbb{E}}  

\newcommand{\Nbb}{\mathbb{N}} 
 
\newcommand{\Pbb}{\mathbb{P}} 
\newcommand{\ud}{{\rm d}}
 
\newcommand{\Prm}{{\rm P}}
\newcommand{\Krm}{{\rm J}}
\newcommand{\Rvrm}{{\rm Rv}}

\newcommand{\Vcal}{\mathcal{V}}

\newcommand{\Scal}{\mathcal{S}}
\newcommand{\Bcal}{\mathcal{B}}
\newcommand{\Ccal}{\mathcal{C}}

\begin{document}

\begin{frontmatter}



\title{Markov chains, AR models, and regular variation}

\author[label1]{Piotr Dyszewski \fnref{label3}} 
 \ead{pdysz@math.uni.wroc.pl}
\fntext[label3]{PD was supported by the National Science Centre, Poland 
(Sonata, grant no.~2020/39/D/ST1/00258).}

\author[label1]{Tamara Mika \corref{cor1}}
 \ead{tamara.fraczek@math.uni.wroc.pl}
\cortext[cor1]{Corresponding author}

 \affiliation[label1]{organization={Mathematical~Institute, University~of~Wroclaw},
             addressline={pl.~Grunwaldzki~2},
             city={Wroclaw},
             postcode={50-384},
             country={Poland}}

\begin{abstract}
	We investigate multivariate regular variation in the context of 
	time-homogeneous Markov chains on general vector spaces and in random coefficient 
	linear models. 
    
    In the first part, we show that the regular variation of the 
	stationary distribution can be derived from that of the innovations, provided that 
	the chain satisfies a certain monotonicity condition with respect to a gauge 
	function. 

	In the second part, we study random linear models with random coefficients 
	defined by an explicit iterative scheme. We prove that the precise structure of 
	the underlying chain affects the form of the associated spectral measure. 
\end{abstract}



\end{frontmatter}



\section{Introduction}

\subsection{The synopsis}
  Regular variation is a fundamental concept in the analysis of heavy-tailed phenomena, 
  which arise in fields such as finance, insurance, and meteorology. 
  The theory was originally developed for univariate functions and real-valued random variables; 
  see, for example, Bingham et al.~\cite{bingham1989regular} and Resnick~\cite{resnick2008extreme}. 
  It was later extended to random vectors~\cite{de1977limit,resnick2007heavy,BASRAK20091055,rootzen2006multivariate}, 
  and this development naturally led to the study of regular variation in function spaces~\cite{Resnick_1986,de2001}.

  Important contributions in this direction include the work of Hult and Lindskog~\cite{hult2006regular}, 
  who introduced $M_0$-convergence for the analysis of measures on metric spaces equipped with scalar multiplication. 
  Building on these ideas, Meinguet and Segers~\cite{meinguet2010regularly} and Segers et al.~\cite{segers2017polar} 
  developed a systematic framework for regularly varying time series in Banach spaces, 
  extending classical results and introducing structural tools such as polar decompositions in 
  star-shaped metric spaces. For a comprehensive introduction to regular 
  variation and its modern developments, we refer to~\cite{kulik2020heavy,mikosch2024extreme}.

  A theoretical study of heavy-tailed phenomena typically requires the analysis of stochastic processes, 
  which are often modeled as Markov chains or as functionals of Markov chains. 
  A useful tool for studying the extremal behavior of Markov chains is the \emph{tail chain}~\cite{planinic2018tail}, 
  which arises under a regularity assumption on the transition kernel~\cite{resnick2013asymptotics}. 
  Heuristically, this condition means that the Markov chain is asymptotically multiplicative. 
  In this setting, one usually works with a stationary version of the chain~\cite{basrak2002regular}.

  In this context, the study of Markov kernels was initiated by Smith~\cite{Smith_1992} in the univariate 
  setting and extended to the multivariate case by Perfekt~\cite{perfekt_1994,Perfekt_1997}. 
  The primary objective of the present article is to analyze how regular variation is propagated by a Markov 
  transition kernel outside the stationary setting, and then to use this analysis to 
  describe the regular variation of the stationary measure. 
  To the best of our knowledge, arguments of this type have so far been developed mainly for 
  real-valued chains~\cite{grey1994regular,Palmowski_Zwart_2007,damek2018iterated,kesten1973random,goldie1991implicit} or 
  for very specific classes of multivariate Markov chains~\cite{AAP1579}.

  This connects with heavy-tailed time-series models where 
  stochastic processes are often represented through infinite series of 
  random vectors~\cite{embrechts2013modelling}. Many models, including stochastic integrals and solutions 
  to random difference equations, admit representations of this type~\cite{mikosch2000supremum}, typically of the form
  \begin{equation*}
    \sum_{j=1}^\infty M_j N_j,
  \end{equation*}
  where $\{N_j\}_{j \in \mathbb{N}}$ is an i.i.d.\ sequence of innovations and 
  $\{M_j\}_{j \in \mathbb{N}}$ is a sequence of random coefficients.

  The Markov-kernel perspective becomes useful when the random coefficients are generated 
  by an underlying Markovian mechanism. In that case, the tail behavior of the coefficient 
  sequence $\{M_j\}_{j\in \mathbb{N}}$, 
  and hence the spectral measure of the resulting series, depends on the structure of the associated Markov chain.

  The heavy-tailed behavior of random-coefficient linear models has been studied 
  extensively~\cite{hult2008tail,davis1995point}. When the innovations are heavy-tailed, 
  large values of the infinite sum are typically caused by a small number of 
  exceptionally large noise variables, whose effect is amplified or attenuated by the random coefficients. 
  In most analyses of such models, the innovations 
  $\{N_j\}_{j \in \mathbb{N}}$ are assumed to be heavy-tailed, whereas the random 
  coefficients $\{M_j\}_{j \in \mathbb{N}}$ are relatively light-tailed.

  Our aim is to relax this assumption by allowing the coefficients themselves to be heavy-tailed. 
  In this setting, the dependence structure of $\{M_j\}_{j \in \mathbb{N}}$ plays a crucial role 
  in the extremal behavior of the series~\cite{mikosch2016large}. 
  Understanding this mechanism is important in applications where joint extremes are central, 
  such as simultaneous crashes in financial markets, catastrophic losses in insurance portfolios, 
  and heavy traffic bursts in telecommunication networks.

\subsection{The contributions and organization of the article}

  The article develops results in three main directions: multivariate regular variation in Section~\ref{sec:mrv}, 
  heavy-tailed Markov chains in Section~\ref{sec:MC}, and random-coefficient linear models in Section~\ref{sec:ar}. 
  We briefly summarize the main contributions of these sections.

  Section~\ref{sec:mrv} recalls the notion of multivariate regular variation in Banach spaces. 
  Its main result, Lemma~\ref{lem:2:equiv}, shows that the corresponding $M_0$-convergence can be 
  verified on a suitable class of test functions. This criterion will be used throughout the paper.

  Section~\ref{sec:MC} applies the results of the preceding section to non-stationary Markov chains with 
  heavy-tailed innovations, in the setting where the transition mechanism is asymptotically homogeneous. 
  In addition to analyzing the associated tail chain, we prove in Theorem~\ref{thm:2:pi} that, 
  for ergodic chains satisfying a suitable monotonicity condition, the stationary distribution is regularly varying.

  Section~\ref{sec:ar} investigates random-coefficient linear models in which both the innovations and 
  the coefficients may be heavy-tailed. We show that, when the coefficients are generated by a Markov chain, 
  the joint regular variation of the coefficients and innovations affects the extremal behavior of the entire system, 
  including the resulting index of regular variation.

\section{Multivariate regular variation}\label{sec:mrv}

  In this section, we recall the notation and the form of multivariate regular variation used throughout the paper. 
  The material is standard; see, for example, \cite{hult2006regular,segers2017polar} for formulations in general 
  metric spaces, and \cite{ledoux1991probability,araujo1980central} for background on probability in Banach spaces. 
  We work in the setting of separable Banach spaces, although several of the statements below extend to 
  star-shaped metric spaces. The main point at which the Banach-space setting differs from the finite-dimensional 
  Euclidean case is that the unit sphere need not be compact; 
  this will be reflected in the choice of test functions used below.

  Throughout the article, let $(\Vcal,|\cdot|)$ be a separable Banach space over $\mathbb{R}$. We denote by
  \begin{equation*}
    \Scal = \{\xbf \in \Vcal : |\xbf|=1\}
  \end{equation*}
  the unit sphere in $\Vcal$, and by
  \begin{equation*}
    \Bcal = \{\xbf \in \Vcal : |\xbf|\le 1\}
  \end{equation*}
  the closed unit ball.

  \begin{df}\label{def:2:RV1}
    A $\Vcal$-valued random element $\Xbf$ is said to be \textit{regularly varying} with index $\alpha>0$ if
    \begin{equation*}
      \Pbb\left[
      \frac{\Xbf}{\|\Xbf\|} \in \cdot , \:
      \frac{\|\Xbf\|}{t} \in \cdot
      \middle| \|\Xbf\|>t
      \right]
      \to
      \Pbb[\Theta \in \cdot]\Pbb[Z \in \cdot],
      \qquad t\to\infty,
    \end{equation*}
    weakly in $\Scal \times (0,+\infty)$, where $\Theta$ is an $\Scal$-valued random element and $Z$ is a 
    Pareto random variable satisfying
    \begin{equation*}
      \Pbb[Z>r]=r^{-\alpha}, \qquad r\ge 1.
    \end{equation*}
    In this case we write $\Xbf \in \Rvrm(\alpha,\Theta)$.
  \end{df}

  In particular, the tail of $\|\Xbf\|$ is regularly varying:
  \begin{equation}\label{urv}
    \frac{\Pbb[\|\Xbf\|>cx]}{\Pbb[\|\Xbf\|>x]} \to c^{-\alpha},
    \qquad x\to\infty,\ c>0.
  \end{equation}
  Equivalently, $\Pbb[\|\Xbf\|>t]=t^{-\alpha}L(t)$ for some slowly varying function $L$. Let $a_t$ be chosen so that
  \begin{equation}\label{eq:2:at}
    t\Pbb[\|\Xbf\|>a_t] \to 1 .
  \end{equation}
  Then $(a_t)_t$ is regularly varying with index $1/\alpha$.

  By \cite[Proposition~3.1]{segers2017polar}, Definition~\ref{def:2:RV1} is equivalent to regular variation in the 
  sense of Hult and Lindskog~\cite{hult2006regular}, formulated in terms of $M_0$-convergence. 
  To state it we recall the corresponding class of test functions. Let
  \begin{equation*}
    \Vcal_0=\Vcal\setminus \{0\}.
  \end{equation*}
  Denote by $M_0(\Vcal)$ the class of Borel measures on $\Vcal_0$ that are finite on $\epsilon\Bcal^c$ 
  for every $\epsilon>0$. Let $\Ccal_0$ be the class of bounded continuous 
  functions $f\colon \Vcal_0\to\mathbb{R}$ that vanish on $\epsilon\Bcal$ for some $\epsilon>0$. 
  For $\mu_n,\mu\in M_0(\Vcal)$, we write
  \begin{equation*}
    \mu_n\to\mu \quad\text{in }M_0(\Vcal)
  \end{equation*}
  if
  \begin{equation*}
    \int f \ud\mu_n\to\int f \ud\mu,
    \qquad f\in\Ccal_0.
  \end{equation*}

  As already mentioned \cite[Proposition~3.1]{segers2017polar}, Definition~\ref{def:2:RV1} is equivalent to
  \begin{equation*}
    t\Pbb[a_t^{-1}\Xbf\in\cdot]\to \mu(\cdot)
    \quad\text{in }M_0(\Vcal),
  \end{equation*}
  where the limit measure $\mu$ is given by
  \begin{equation}\label{eq:2:muTheta}
    \mu(A)=\int_0^\infty \Pbb[s\Theta\in A]\alpha s^{-\alpha-1}\ud s .
  \end{equation}
  We note that the measure $\mu$ satisfies
  \begin{equation}\label{eq:2:norm}
    \mu(\Vcal\setminus\Bcal)=1
  \end{equation}
  and the homogeneity relation
  \begin{equation}\label{eq:2:hom}
    \mu(sA)=s^{-\alpha}\mu(A),
    \qquad s>0.
  \end{equation}

  For later use, we also introduce a smaller class of test functions. 
  Let $\Ccal_0^u$ be the class of bounded, nonnegative, continuous functions $f:\Vcal\to\mathbb{R}$ 
  such that $0\notin\operatorname{supp}(f)$ and whose restrictions to $M\Bcal$ are uniformly continuous for every $M>0$. 
  This additional uniform-continuity requirement is harmless in finite dimensions, 
  where bounded closed balls are compact, but it is useful in the Banach-space setting, where such compactness is not available.

  The next lemma shows that these test functions are sufficient to characterize the above convergence.

  \begin{lem}\label{lem:2:equiv}
    Let $\Xbf$ be a $\Vcal$-valued random vector. The following statements are equivalent:
    \begin{itemize}
      \item[(i)] $\Xbf \in \Rvrm(\alpha,\Theta)$.
	
      \item[(ii)] Let $\mu$ be defined by~\eqref{eq:2:muTheta}. Then, for any regularly varying
		    function $a$ satisfying~\eqref{eq:2:at},
        \begin{equation*}
    		  t\,\Pbb[a_t^{-1}\Xbf \in \cdot] \to \mu(\cdot)
    		  \qquad \text{in } M_0(\Vcal).
        \end{equation*}

		  \item[(iii)] For every $f\in\Ccal_0^u$ and any regularly varying function $a$
		    satisfying~\eqref{eq:2:at},
		    \begin{equation}\label{eq:2:claim}
    		  t\,\Ebb\!\left[f\!\left(a_t^{-1}\Xbf\right)\right]
    		  \to \int f(\xbf)\,\mu(\ud\xbf),
		    \end{equation}
		    where $\mu$ is defined by~\eqref{eq:2:muTheta}.
    \end{itemize}
  \end{lem}

  \begin{proof}
    The equivalence of (i) and (ii) follows from~\cite[Proposition~3.1]{segers2017polar}. 
    It remains to prove the equivalence of (i) and (iii).

    \medskip
    \noindent
    (i) $\Rightarrow$ (iii).
    Assume that $\Xbf \in \Rvrm(\alpha,\Theta)$ and fix $f\in\Ccal_0^u$. 
    Then there exists $\epsilon>0$ such that $f(\xbf)=0$ whenever $\|\xbf\|\le\epsilon$. 
    Define $g\colon (0,\infty)\times\Scal\to[0,\infty)$ by
    \[
		g(r,\theta)=f(\epsilon r\theta).
    \]
    The function $g$ is bounded and continuous. Hence, by Definition~\ref{def:2:RV1},
    \[
		\Ebb\!\left[f\!\left(\epsilon t^{-1}\Xbf\right)\,\middle|\,\|\Xbf\|>t\right]
		\to \int_1^\infty \Ebb\!\left[f(\epsilon r\Theta)\right]\alpha r^{-\alpha-1}\,\ud r.
    \]
    We now prove~\eqref{eq:2:claim}. Since $f(\xbf)=0$ for $\|\xbf\|\le \epsilon$, we have
    \begin{align*}
	   \Ebb\!\left[f\!\left(a_t^{-1}\Xbf\right)\right]
	   & = \Ebb\!\left[f\!\left(a_t^{-1}\Xbf\right)\mathbf{1}_{\{\|\Xbf\|>\epsilon a_t\}}\right]\\
	   & = \Pbb[\|\Xbf\|>\epsilon a_t]\,
	   \Ebb\!\left[f\!\left(a_t^{-1}\Xbf\right)\,\middle|\,\|\Xbf\|>\epsilon a_t\right].
    \end{align*}
    By~\eqref{eq:2:at} and~\eqref{urv},
    \[
	   t\,\Pbb[\|\Xbf\|>\epsilon a_t]
	   \to \epsilon^{-\alpha}.
    \]
    Therefore
    \[
	   t\Ebb\!\left[f\!\left( a_t^{-1}\Xbf\right)\right]
	   \to \epsilon^{-\alpha} \int_1^\infty \Ebb\!\left[f(\epsilon r\Theta)\right]\alpha r^{-\alpha-1}\,\ud r.
    \]
    The change of variables $s=\epsilon r$ gives~\eqref{eq:2:claim}.

    \medskip
    \noindent
    (iii) $\Rightarrow$ (i).
    Assume that~\eqref{eq:2:claim} holds for every $f\in\Ccal_0^u$. 
    Let $\Ccal'$ denote the class of all bounded functions
    $g\colon (0,\infty)\times\Scal\to[0,\infty)$ that are continuous with respect to the product topology,
    uniformly continuous on $[1,M]\times\Scal$ for every $M>0$, 
    and supported in $[1,\infty)\times\Scal$.

    For each $g\in\Ccal'$, define $f(0)=0$ and for $x \neq 0$,
    $f(\xbf)=g(\|\xbf\|,\xbf/\|\xbf\|)$.
    Then $f\in\Ccal_0^u$. Hence,
    \[
		\Ebb\!\left[g\!\left(t^{-1}\|\Xbf\|,\Xbf/\|\Xbf\|\right)
		\,\middle|\,\|\Xbf\|>t\right]
		\to \Ebb[g(Z,\Theta)],
    \]
    where $Z$ has a Pareto distribution with $\Pbb[Z>r]=r^{-\alpha}$ for $r\ge1$. 
    Since $\Ccal'$ contains all uniformly continuous functions supported in 
    $[1,\infty)\times\Scal$, this yields the convergence required in Definition~\ref{def:2:RV1}.
  \end{proof}

  The previous lemma motivates the following equivalent definition of 
  multivariate regular variation.

  \begin{df}\label{df:4:RV2}
	  Let $\mu$ be a measure on $\Vcal_0$ such that \eqref{eq:2:norm} and~\eqref{eq:2:hom} hold. 
	  A $\Vcal$-valued random vector $\Xbf$ is said to be \emph{regularly varying} 
	  with index $\alpha>0$ and limit measure $\mu$ if there exists a function 
	  $a:(0,\infty)\to(0,\infty)$, regularly varying with index $1/\alpha$, such that
	  \[
		  t\,\Pbb[a_t^{-1}\Xbf\in\cdot]\to\mu(\cdot)
		  \qquad \text{in } M_0(\Vcal).
	  \]
	  In this case we write $\Xbf\in\Rvrm(\alpha,\mu)$.
  \end{df}

  Lemma~\ref{lem:2:equiv} shows that Definitions~\ref{def:2:RV1} and~\ref{df:4:RV2}
  are indeed equivalent.

  It is worth noting that, when $\Vcal$ is finite-dimensional, multivariate regular
  variation can also be formulated in terms of vague convergence. In that case, this amounts to requiring that
  \begin{equation*}
    \lim_{t\to\infty} t\Ebb\left[h\left(a_t^{-1}\Xbf\right)\right]
    = \int h(\xbf)\mu(\ud\xbf),
  \end{equation*}
  for every continuous, nonnegative function $h$ with compact support in $\Vcal_0$.
  However, this approach does not extend naturally to infinite-dimensional spaces.
  The difficulty is that vague convergence relies on local compactness, which
  infinite-dimensional Banach spaces lack. This limitation becomes apparent in the
  proof of the implication (iii) $\Rightarrow$ (i) in Lemma~\ref{lem:2:equiv}, where the
  constructed function $f$ need not have compact support.

\section{Markov chains}\label{sec:MC}

  We now introduce the framework used in our results 
  on the regular variation of $\Vcal$-valued Markov chains. 
  Throughout, let $\Sigma$ denote the Borel $\sigma$-algebra 
  on $\Vcal$, i.e., the $\sigma$-algebra generated by all open subsets of $\Vcal$. 

\subsection{Heavy-tailed transition kernel}

  Recall that a function $\Prm \colon \Vcal \times \Sigma\to [0,1]$ 
  is called a transition function if, 
  for every $\xbf \in \Vcal$, the map $A \mapsto \Prm(\xbf, A)$ is a 
  probability measure on $(\Vcal,\Sigma)$, 
  and for every $A \in \Sigma$, the map $\xbf \mapsto \Prm(\xbf, A)$ is Borel measurable.
  We write $\sigma(X_1, X_2, \dots)$ for the $\sigma$-algebra generated by $X_1, X_2, \dots$.

  Consider a family 
  $\{\Xbf_n^{\xbf} : n \in \mathbb{N},\, \xbf \in \Vcal\}$ 
  of time-homogeneous Markov chains on $\Vcal$ 
  with transition function $\Prm$ with the feller property.
  By this we mean a family of $\Vcal$-valued random elements 
  defined on a common probability space and satisfying the following conditions:
  \begin{enumerate}
    \item For every $\xbf \in \Vcal$, $\Pbb[\Xbf_0^{\xbf} = \xbf] = 1$.

    \item For every $n \in \mathbb{N}$ and every $\xbf \in \Vcal$, 
      if $\mathcal{F}_n^{\xbf} = \sigma(\Xbf_j^{\xbf} : 0\leq j \leq n)$,
      then the regular conditional distribution of $\Xbf_{n+1}^{\xbf}$ 
      given $\mathcal{F}_n^{\xbf}$ is given by
      \begin{equation}
   		  \Pbb \left[\Xbf_{n+1}^{\xbf} \in \ud \ybf \;\middle|\; \mathcal{F}_n^{\xbf} \right] 
   		  = \Prm\!\left(\Xbf_n^{\xbf}, \ud \ybf\right).
      \end{equation}
      \item for any bounded and continuous functon $f \colon \Vcal \to \mathbb{R}$,
        the function
        \begin{equation*}
            x \mapsto \int_\Vcal f(y) \Prm(x, \mathrm{d}y)
        \end{equation*}
        is continuous.
  \end{enumerate}

  Our goal is to analyse the effect of heavy-tailed innovations on the extremal 
  behaviour of the Markov chain $\{\Xbf_n^{\xbf}\}_{n \in \mathbb{N}}$, 
  especially in settings where the one-big-jump principle applies.

  \begin{as}\label{as:2:first}
    There exist $\xbf_0 \in \Vcal$ and a positive, non-decreasing function 
    $a \colon (0,\infty) \to (0,\infty)$, 
    regularly varying with index $1/\alpha$, such that 
    $\mathbb{E}\left[ \left\|\Xbf_1^{\xbf_0}\right\|^\alpha \right]<\infty$,
    \begin{equation}\label{eq:2:banana}  
    	\lim_{t \to \infty} t \Pbb[\|\Xbf_1^{\xbf_0}\| > a_t] = 1,  
    \end{equation}  
    and the following four conditions hold:
    \begin{itemize}
      \item [$(i)$] \textbf{Innovations:} There exists a family 
    	  $\{\nu_\bfx\}_{\bfx \in \Vcal}$ of measures in $M_0(\Vcal)$ such that, 
    		for every $\bfx \in \Vcal$, $\nu_x(\Bcal^c)>0$ and
    		\begin{equation} \label{eq:4:gass1}  
    		  \lim_{t \to \infty} t \Prm(\bfx, a_t \cdot) = \nu_\bfx(\cdot)  
    		\end{equation}  
    		in $M_0(\Vcal)$.  

      \item [$(ii)$] \textbf{Asymptotic stability:} There exists a transition kernel
    		$\Krm$ such that, for every $\bfx \in \Vcal$,  
    		\begin{equation} \label{eq:4:gass2}  
    		  \lim_{t \to \infty} \Prm(t \bfx, t \cdot) = \Krm(\bfx,\cdot)  
    		\end{equation}  
    		weakly, with convergence uniform on bounded sets. 
    		More precisely, for every $R > 0$ and every $f \in \Ccal_0^u$,  
    		\begin{equation}\label{eq:3:conv}
    		  \lim_{t \to \infty} \sup_{\|\bfx\| \leq R}  
    			\left| \int f(\bfy) \Prm(t \bfx, t \ud \bfy) 
    			- \int f(\bfy) \Krm(\bfx,\ud \bfy) \right| = 0.  
    		\end{equation}  

      \item [$(iii)$] \textbf{Multiplicative domination:} There exists a constant $C>0$ such that
        \begin{equation*}
          \Prm(x, t\mathcal{B}^c) \leq C \Prm(x_0, \|x\|^{-1}t \mathcal{B}^c)
        \end{equation*}
        for all $x \in \mathcal{V}_0$ with
        $\|x\| >C$ and all $t>0$.

    		\item [$(iv)$] \textbf{Single jumps:} For independent random variables
          $X_1$ $X_2$ distributed as $\|X_1^{x_0}\|$ we have  
    			\begin{equation} \label{gass3}  
    			  \lim_{M \to \infty} \limsup_{t \to \infty} t  
    				\Pbb\big[X_1X_2 \geq a_t, X_1 \in (M, a_t/M)\big] = 0.  
    			\end{equation}  
	  \end{itemize}
  \end{as}

  The first condition~\eqref{eq:4:gass1} in Assumption~\ref{as:2:first} implies 
  that the innovations in the Markov 
  chain under consideration are heavy-tailed. In particular for any $x \in \mathcal{V}$
  the random vector $\Xbf_1^x$ is $\Rvrm(\alpha, \mu)$ with 
  \begin{equation*}
    \mu(\cdot) = \frac{1}{\nu_x\left(\mathcal{B}^c\right)} \nu_x(\cdot).
  \end{equation*}
  In particular the measure $\nu_x$ is homogeneous in the sense of~\eqref{eq:2:hom}.
  Note that then by the third condition
  \begin{equation}\label{eq:3:defR} 
    R(x) = \lim_{t \to \infty} \frac{\Prm(x, t \Bcal^c)}{\Prm(x_0, t \Bcal^c)}
	= \nu_x ( \Bcal^c) \leq C \|x\|^\alpha <\infty
  \end{equation}
  for all $x \in \mathcal{V}_0$ with $\|x\| >C$
  and therefore the tail probabilities of all $X_1^x$'s are comparable.

  The final condition~\eqref{gass3} 
  ensures that the sequence $\{X_n^x\}_{n \in \mathbb{N}}$ 
  adheres to the principle of a single big jump and is similar in nature to the one used for products of random matrices. It implies that, under the notation of Assumption~\ref{as:2:first},
  \begin{equation*}
      \mathbb{P}[X_1X_2>t] \sim 2 \mathbb{E}[X_1^\alpha] \mathbb{P}[X_1>t].
  \end{equation*}
	
  The second condition~\eqref{eq:4:gass2} quantifies the 
  impact of large jumps on the future evolution of the process, 
  aligning with the condition proposed by~\cite{resnick2013asymptotics} in 
  the context of the real line. 
  Notably, this framework reveals a Markovian structure in the extremes of $\Xbf^\xbf$.
  To describe this, let $\{\Zbf_n^\xbf \: : \: n\in \mathbb{N}, \: x \in \Vcal\}$
  be a Markov chain with transition function $\Krm(\xbf, \ud \ybf)$.
  Note that $\{\Zbf_n^\xbf\}_{ n \in \Nbb}$ is homogeneous in $\xbf$. 
  Indeed, for any $c>0$ and any $\xbf \in \Vcal$ by definition 
  one must have that 
  \begin{equation*}
    \Krm(c \xbf, \ud \ybf) = \Krm\left( \xbf, c^{-1}\ud \ybf \right).
  \end{equation*}
  This by Markov property implies that
  the processes $\{\Zbf_n^{c\xbf}\}_{n \in \mathbb{N}}$ and 
  $\{c\Zbf_n^{\xbf}\}_{n \in \mathbb{N}}$ have the same law. 
  As we will see the chain
  $\{\Zbf_n^{\xbf}\}_{n \in \mathbb{N}}$ arises in the description of the extremes
  of $\{\Xbf_n^{\xbf}\}_{n \in \mathbb{N}}$.

  \begin{example}
    Let $\Vcal = \mathbb{R}^d$ and
    consider a Markov chain given via so-called random difference equation
    \begin{equation*}
      \Xbf_{n+1}^\xbf = A_{n+1} \Xbf_n^{\xbf} + B_{n+1},
    \end{equation*}
    where for each $n \in \mathbb{N}$, 
    $A_{n+1}$ is $d \times d$ random matrix and $B_{n+1}$ is a 
    $d$-dimensional random vector such that $\{ (A_n, B_n)\}_{n \in \mathbb{N}}$
    are iid. Assume that the pair $(A_1,B_1)$ is $\Rvrm(\alpha, \mu)$ 
    in the space 
    $M_{d\times d}(\mathbb{R}) \times \mathbb{R}^d$, 
    where $M_{d\times d}(\mathbb{R})$ denotes
    the space of $d\times d$ matrices with real entries. 
    Equip $M_{d\times d}(\mathbb{R}) \times \mathbb{R}^d$ with the norm
    $\|(a,b)\|=\|a\| + \|b\|$, where $\|a\|$ is the spectral norm of $a$ and
    $\|b\|$ is the euclidean norm of $\|b\|$. This choice of norm will allow us to get the multiplicative domination condition for this Markov chain.
    In what follows we assume that
    \begin{equation*}
        \mathbb{E}[\|(A_1,B_1)\|^\alpha]<\infty.
    \end{equation*}
    From the condition \eqref{eq:2:norm} we have
    \begin{equation*}
      \mu \left( (a,b) \in M_{d \times d}(\mathbb{R}) \times \mathbb{R}^d \: : 
      \: \|(a,b)\|\geq 1\right)=1.
    \end{equation*}
    There has to exist $x_0 \in \mathbb{R}^d$ such that
    \begin{equation*}
      L_{x_0}=\mu \left( (a,b) \in M_{d \times d}(\mathbb{R}) \times \mathbb{R}^d \: : 
      \: \|ax_0+b\|\geq 1\right)>0.
    \end{equation*}
    Then the random vector $X_1^{x_0}=A_1x_0+B_1$ is regularly varying 
    and there exists $a \colon \mathbb{R} \to \mathbb{R}$ satisfying~\eqref{eq:2:banana}. 
    It follows that~\eqref{eq:4:gass1} holds with
    \begin{equation*}
      \nu_x(\cdot) = L_{x_0}^{-1} 
      \mu \left( (a,b) \in M_{d \times d}(\mathbb{R}) \times \mathbb{R}^d \: : 
      \: ax+b \in \cdot\right).
    \end{equation*}
    Since 
    \begin{equation*}
      \Prm(tx, t \cdot) = \mathbb{P} \left[ A_1x + t^{-1}B_1 \in \cdot  \right]
      \to \Krm(x, \cdot) = \mathbb{P} \left[ A_1x \in \cdot  \right]
    \end{equation*}
    weakly, with the convergence being uniform in $x$ in the sense 
    of \eqref{eq:3:conv}. Therefore the Markov chain 
    $\{Z_n^x \: : \: x \in \mathbb{R}^d, \: n \in \mathbb{N}\}$ is a multiplicative 
    random walk
    \begin{equation*}
      Z_n^x =  A_n A_{n-1} \cdots A_1 x.
    \end{equation*}
    In summary, the first three conditions~\eqref{eq:2:banana}-\eqref{eq:4:gass2} 
    of Assumption~\ref{as:2:first} are satisfied if $(A_1, B_1)$ is regularly 
    varying in $M_{d\times d}(\mathbb{R}) \times \mathbb{R}^d$.
    The last condition~\eqref{gass3} is satisfied if additionally one assumes that
    \begin{equation*}
      \lim_{M \to \infty} \limsup_{t \to \infty} t 
      \mathbb{P}\left[ D_1 D_2>a_t, \: D_1 \in (M, a_t/M) \right] = 0,
    \end{equation*}
    where $D_i = \|(A_i,B_i)\|$. Finally the multiplicative domination condition 
    follows since for $\|x\|\geq 1$ we have
    \begin{align*}
      \mathbb{P}\left[\|Ax+B\| >t\right] 
        & \leq \mathbb{P}\left[ \|(A,B)\| \cdot (\|x\|+1) >t \right] \\
        & \leq C \mathbb{P}\left[ \|Ax_0+B\| \cdot (\|x\|+1) >t \right] \\
        & \leq C \mathbb{P}\left[ \|Ax_0+B\| \cdot \|x\| >t \right].
    \end{align*}
    The first inequality is a consequence of the tail equivalence of $\|(A, B)\|$ and 
    $\|Ax_0+B\|$ and the second inequality 
    holds by possibly increasing the value of $C$.
  \end{example}

\begin{example}
Let $\Vcal=C([0,1])$, $\|f\|_\infty=\sup_{u\in[0,1]}|f(u)|$,
and let $\mathcal A\colon\Vcal\to\Vcal$ be a bounded linear operator
with $|\mathcal A|<1$. Let
$(\varepsilon_n)_{n\geq1}$ be an i.i.d.\ sequence of regularly varying
random elements of $C([0,1])$, and consider the Markov chain
\begin{equation*}
\Xbf_{n+1}^{\xbf}
=
\mathcal A\Xbf_n^{\xbf}+\varepsilon_{n+1},
\qquad
\Xbf_0^{\xbf}=\xbf.
\end{equation*}
This is a functional autoregressive process of order one; see
\cite{bosq2000linear}. Such processes are used to model time series
whose individual observations are curves rather than finite-dimensional
vectors.
The transition kernel is
$\Prm(\xbf,\cdot)
=
\mathcal L(\mathcal A\xbf+\varepsilon_1)$,
while its asymptotic transition kernel is
$\Krm(\xbf,\cdot)=\delta_{\mathcal A\xbf}$.
Hence the corresponding tail chain is simply
$\Zbf_n^\xbf=\mathcal A^n\xbf$.
\end{example}

\subsection{Regular variation in Markov chain}

  We now move to the analysis of regular variation of $X_n^x$'s. We begin with a lemma that asserts
  that the regular variation property is preserved by $\Prm$.

  \begin{lem}\label{lem:3:key}
    Let Assumption~\ref{as:2:first} be in force.
    Assume that a $\Vcal$-valued random vector 
    $\Xbf$ is in $\Rvrm(\alpha, \mu_\Xbf)$ for some measure $\mu_\Xbf$. 
    Denote  
    \begin{equation*}
      \lim_{t \to \infty} t \Pbb[\|\Xbf\| > a_t] = L_\Xbf.  
    \end{equation*}  
    Consider a $\Vcal$-valued random vector $\Ybf$ with conditional distribution
    \begin{equation*}  
      \Pbb\left[ \left. \Ybf \in \cdot \right| \Xbf \right] = \Prm(\Xbf, \cdot).
    \end{equation*}  
    then $\Ybf$ is in $\Rvrm(\alpha, \mu_\Ybf)$, where  
    \begin{equation*}  
      L_\Ybf \mu_\Ybf(\cdot) = 
      \Ebb[\nu_{\Xbf}(\cdot)] + L_\Xbf \int \Krm(\bfy,\cdot) \mu_\Xbf(\ud \bfy),  
    \end{equation*}  
    and  
    \begin{equation*}  
      L_\Ybf = \lim_{t \to \infty} t \Pbb[\|\Ybf\| > a_t] =  
      \Ebb[\nu_{\Xbf}(\Bcal^c)] + L_\Xbf \int \Krm(\bfy, \mathcal{B}^c) \mu_\Xbf(\ud \bfy).  
    \end{equation*}
  \end{lem}

  \begin{proof}
    For \( M \geq 1 \), let \(\varphi_M \colon \Vcal \to [0,1]\) 
    be any continuous function such that  
    \begin{equation*}  
      \varphi_M(\bfy) =  
      \begin{cases}  
        1, & \|\bfy\| \leq M, \\  
        0, & \|\bfy\| \geq 2M.  
      \end{cases}  
    \end{equation*}  
    Similarly, let \(\psi_M \colon \Vcal \to [0,1]\) 
    be any continuous function such that  
    \begin{equation*}  
      \psi_M(\bfy) =  
      \begin{cases}  
        1, & \|\bfy\| \geq 1/M, \\  
        0, & \|\bfy\| \leq 1/(2M).  
      \end{cases}  
    \end{equation*}  
    To prove that \(\Ybf\) is regularly varying, 
    take \( f \in \Ccal_0^u \) and write  
    \begin{multline*}  
      t \Ebb \left[f\left(a_t^{-1} \Ybf\right)\right] =  
      t \Ebb\left[f\left(a_t^{-1} \Ybf\right) \varphi_M\left(\Xbf\right)\right] + \\  
      t \Ebb\left[f\left(a_t^{-1} \Ybf\right) \left(1 - \varphi_M\left(\Xbf\right)\right)  
      \left(1 - \psi_M\left(a_t^{-1} \Xbf\right)\right)\right] + \\  
      t \Ebb\left[f\left(a_t^{-1} \Ybf\right) \psi_M\left(a_t^{-1} \Xbf\right)\right] 
      =: D_1 + D_2 + D_3,  
    \end{multline*}  
    for sufficiently large $t$.  

    The first term can be reformulated as  
    \begin{equation*}  
      D_1 = \Ebb \left[t \int f(\bfy) \Prm(\Xbf, a_t \ud \bfy) \varphi_M(\Xbf)\right].  
    \end{equation*}  
    Using~\eqref{eq:4:gass1} 
    and the dominated convergence theorem, we obtain  
    \begin{equation*}  
      \lim_{t \to \infty} D_1 
      = \Ebb \left[\int f(\bfy) \nu_{\Xbf}(\ud \bfy) \varphi_M(\Xbf)\right].  
    \end{equation*}  
    Applying the monotone convergence theorem once again, we conclude that  
    \begin{equation*}  
      \lim_{M \to \infty} \lim_{t \to \infty} D_1 
      = \Ebb \left[\int f(\bfy) \nu_{\Xbf}(\ud \bfy)\right].  
    \end{equation*}  

    The second part of the decomposition, $D_2$, 
    is negligible as $M \to \infty$ and $t \to \infty$. 
    This follows from the fact that there exists $\epsilon > 0$ such that  
    $\{\|\xbf\| < \epsilon\} \cap \operatorname{supp}(f) = \emptyset$.  
    Since $f$ is non-negative and bounded above by some constant $c > 0$, we have  
    \begin{multline*}  
      D_2 \leq t \Ebb \left[f\left(a_t^{-1} \Ybf\right) \1_{\{M < \|\Xbf\| < a_t / M\}}\right]\\
      \leq c t \Pbb\left[a_t^{-1} \|\Ybf\| > \epsilon, M < \|\Xbf\| < a_t / M\right]\\
      \leq C c t \Ebb \left[ \Prm\left(x_0, \epsilon a_t \|\Xbf\|^{-1}\mathcal{B}^c \right) \1_{\{M < \|\Xbf\| < a_t / M\}}\right]\\
      = C c t \Pbb\left[ \|X_1^{x_0}\| \cdot \|X\| > \epsilon a_t, M < \|\Xbf\| < a_t / M\right]
    \end{multline*}
    Note that by 
    Using \eqref{gass3}, it follows that  
    \begin{equation*}  
      \lim_{M \to \infty} \limsup_{t \to \infty} D_2 = 0.  
    \end{equation*}  

    The third part can be rewritten as  
    \begin{equation*}  
      D_3 = \Ebb\left[t \int f(\bfy) \Prm(\Xbf, a_t \ud \bfy) 
      \psi_M\left(a_t^{-1} \Xbf \right)\right].  
    \end{equation*}  
    Define  
    \begin{equation*}  
      G_t(\bfx) = \int f(\bfy) \Prm(a_t \bfx, a_t \ud \bfy).  
    \end{equation*}  
    Using \eqref{eq:4:gass2}, we have  
    \begin{equation*}  
      G_\infty(\bfx) := \lim_{t \to \infty} G_t(\bfx) = \int f(\bfy) \Krm(\bfx,\ud \bfy),  
    \end{equation*}  
    where the convergence is in the sense of \eqref{eq:3:conv}.
    With this notation, we can express \(D_3\) as  
    \begin{multline*}  
      D_3 =  
      \Ebb\left[t G_t\left(a_t^{-1} \Xbf\right) \left(1 - \varphi_M(a_t^{-1} \Xbf)\right) 
      \psi_M(a_t^{-1} \Xbf)\right] \\  
      + 
      \Ebb\left[t \left(G_t\left(a_t^{-1} \Xbf\right) - G_\infty\left(a_t^{-1} \Xbf\right)\right)  
      \varphi_M\left(a_t^{-1} \Xbf\right) \psi_M\left(a_t^{-1} \Xbf\right)\right] \\  
      + 
      \Ebb\left[t G_\infty\left(a_t^{-1} \Xbf\right) 
      \varphi_M\left(a_t^{-1} \Xbf\right) \psi_M\left(a_t^{-1} \Xbf\right)\right] 
      =: E_1 + E_2 + E_3.  
    \end{multline*}  
    The contribution of \( E_1 \) in \( D_3 \) is negligible. 
    Indeed, since \(\Prm(x, \cdot)\) are probability measures and \( f \) is bounded, we have  
    \begin{equation*}  
      \sup_{t > 0} \sup_{\bfx \in \Vcal} G_t(\bfx) < \infty.  
    \end{equation*}  
    In view of \(\Pbb[\|\Xbf\| \ge M a_t] \sim L_\Xbf M^{-\alpha} t^{-1}\), the last expression yields  
    \begin{equation*}  
      \lim_{M \to \infty} \limsup_{t \to \infty} E_1 = 0.  
    \end{equation*}  
    The term $E_2$ can be handled similarly. 
    We have
    \begin{equation*}
      E_2 \le \sup_{\|\xbf\| \le 2M} \left|G_t(\xbf) - G_\infty (\xbf)\right| \:  
      t \: \P[\|\Xbf\| \ge a_t/ M].
    \end{equation*}
    Using the fact that for all $M>0$
    \begin{equation*}  
      \lim_{t \to \infty} \sup_{\|\bfx\| \leq M} \left| G_t(\bfx) - G_\infty(\bfx) \right| = 0 
    \end{equation*}  
    we have 
    \begin{equation*}
      \lim_{t\to \infty} E_2 \le L_{\Xbf} M^{\alpha} 
      \lim_{t\to \infty} \sup_{\|\xbf\| \le 2M} \left|G_t(\xbf) - G_\infty (\xbf)\right| =0.
    \end{equation*}
    It remains to analyse \( E_3 \). By the regular variation of \(\Xbf\), we have  
    \begin{equation*}  
      \lim_{t \to \infty} E_3 = 
      L_\Xbf \int G_\infty(\bfy) \varphi_M\left(\bfy\right) 
      \psi_M\left(\bfy\right) \mu_\Xbf(\ud \bfy).  
    \end{equation*}  
    As \( M \to \infty \), the right-hand side converges to  
    \begin{equation*}  
      L_\Xbf \int G_\infty(\bfy) \mu_\Xbf(\ud \bfy) = 
      L_\Xbf \int \int f(\bfx) \Krm(\bfy,\ud \bfx) \mu_\Xbf(\ud \bfy).  
    \end{equation*}  
    This completes the proof.
  \end{proof}

  Our first result asserts regular variation of 
  $\Xbf_n^\xbf$'s. We note that in contrast to the 
  literature~\cite{BASRAK20091055,janssen2014markov} we study the extremes of
  $\Xbf_n^\xbf$ conditioned on $\|\Xbf_n^{\xbf}\|$ being large rather than $\|\Xbf_1^\xbf\|$.

  \begin{prop}\label{prop:2:xn}
    Let the Assumption~\ref{as:2:first} be in force. For any 
    $n \in \Nbb$ and any $\xbf \in \Vcal$ the law of $\Xbf_n^\xbf$ 
    is regularly varying. More precisely 
    $\Xbf_n^\xbf\in \Rvrm(\alpha, \mu^{\xbf}_n)$, 
    where
    \begin{equation*}
		L_n \mu^\xbf_n (\cdot)
		= \sum_{j=0}^{n-1} \Ebb \left[ \int \Pbb[ \Zbf_{n-j-1}^\ybf \in \cdot] 
		\nu_{\Xbf^\xbf_j} (\ud \ybf) \right].  
    \end{equation*}
    Here $L_n>0$ is a normalising constant such that $\mu_n(\Bcal^c)=1$. 
  \end{prop}

  \begin{proof}
    Fix $\xbf \in \Vcal$. Then, by~\eqref{eq:4:gass1}, $\Xbf_1^\xbf$ 
    is regularly varying since, 
    \begin{equation*}  
      \lim_{t \to \infty} t \Pbb \left[ a_t^{-1} \Xbf_1^\xbf \in \cdot \right] 
      = \nu_\xbf(\cdot)  
    \end{equation*}  
    in $M_0(\Vcal)$. Thus,  
    \begin{equation*}  
      L_1 = \lim_{t \to \infty} t \Pbb \left[ \|\Xbf_1^\xbf\| \geq a_t \right] 
      = \nu_\xbf(\Bcal^c).  
    \end{equation*} 
    Put $L_1\mu_1^\xbf=\nu_\xbf$. 
    For \( n \in \mathbb{N} \), \( n > 1 \), we apply Lemma~\ref{lem:3:key} 
    to obtain by induction  
    \begin{equation*}  
      \lim_{t \to \infty} t \Pbb \left[ a_t^{-1} \Xbf_n^\xbf \in \cdot \right] 
      = L_n\mu^\xbf_n(\cdot)  
    \end{equation*}  
    in $M_0(\Vcal)$, where  
    \begin{equation*}  
      L_n \mu^\xbf_n(\cdot) =  
      \Ebb[\nu_{\Xbf^\xbf_{n-1}}(\cdot)] + L_{n-1} \int \Krm(\bfy,\cdot) 
      \mu^\xbf_{n-1} (\ud \bfy),  
    \end{equation*}  
    with  
    \begin{equation*}  
      L_n = \lim_{t \to \infty} t \Pbb[\|\Xbf_n^\xbf\| > a_t] = \\  
       \Ebb[\nu_{\Xbf^\xbf_{n-1}}(\Bcal^c)] +  
      L_{n-1} \int \Krm(\bfy, \{\|\bfx\| > 1\}) \mu^\xbf_{n-1} (\ud \bfy).  
    \end{equation*}  
    The recursive relation for \(\{L_n \mu^\xbf_n\}_{n \in \mathbb{N}}\) 
    can be solved as  
    \[
      L_n \mu^\xbf_n (\cdot)
      = \sum_{j=0}^{n-1} \Ebb_\xbf \left[ \int \Pbb[ \Zbf_{n-j-1}^\ybf \in \cdot] 
      \nu_{\Xbf^\xbf_j} (\ud \ybf) \right] 
    \]  
    as claimed.
  \end{proof}

\subsection{Regular variation of the stationary distribution}

  Our second result concerns the regular variation of the stationary distribution 
  of the process $\{\Xbf_n^\xbf\}_{n\in \mathbb{N}}$, denoted by the probability measure 
  $\pi$ on $\Vcal$, which satisfies
  \begin{equation*}
    \pi(\cdot) = \int_{\Vcal} \Prm(\xbf, \cdot)\, \pi(\ud \xbf).
  \end{equation*}
  Moreover, we will need that for every $\xbf \in \Vcal$, 
  the distributions of $\{\Xbf_n^{\xbf}\}_{n \in \mathbb{N}}$ converge weakly to $\pi$ 
  as $n \to \infty$. 
  To ensure this convergence, we impose the following classical
  conditions~\cite{harris1956existence, hairer2011yet}.

  \begin{as}\label{as:2:mix}
    There exists a function $V \colon \Vcal \to [0, +\infty)$ 
    satisfying the following two conditions:  
    \begin{itemize}  
      \item [$(i)$] There exist constants $\kappa > 0$ and $\gamma \in (0,1)$ such that  
        \begin{equation*}  
          \Ebb[V(\Xbf_1^\xbf)] \leq \gamma V(\xbf) + \kappa  
        \end{equation*}  
        for all $\xbf \in \Vcal$.  
      \item [$(ii)$] There exist a constant $\epsilon>0$
        and a probability measure $\nu$ on $\Vcal$ such that
        \begin{equation}\label{eq:2:newbound}  
          \Prm(\xbf, \cdot) \geq \epsilon\nu(\cdot)
        \end{equation}  
        for all $\xbf \in \mathcal{C}$, where 
        $\mathcal{C} = \{ \xbf \in \Vcal \: : \: V(\xbf) \leq R \}$ 
        for some $R> 2\kappa/(1-\gamma)$.
    \end{itemize}
  \end{as}

  Assumption~\ref{as:2:mix} guarantees the existence of a unique stationary 
  distribution $\pi$ and geometric convergence of the distribution of 
  $\Xbf_n^\xbf$ to $\pi$. To make this precise, we introduce the $V$-weighted 
  norm for functions $\phi \colon \Vcal \to \mathbb{R}$, defined by
  \begin{equation*}  
    \|\phi\|_V = \sup_{\xbf \in \Vcal} \frac{|\phi(\xbf)|}{1+V(\xbf)}.  
  \end{equation*}  
  We can now state the following result.

  \begin{prop}\label{prop:2:mix}  
    Suppose that Assumption~\ref{as:2:mix} holds. 
    Then the transition kernel $\Prm$ admits a unique stationary distribution $\pi$. 
    Moreover, there exist constants $\rho \in (0,1)$ and $C_{\pi} > 0$ such that  
    \begin{equation*}  
      \left\| \Ebb \left[\phi\left(\Xbf_n^{(\cdot)}\right)\right] - \int \phi \: \ud \pi 
      \right\|_V \leq  
      C_{\pi} \rho^n \left\| \phi(\cdot) - \int \phi \: \ud \pi \right\|_V  
    \end{equation*}  
    for every $n \in \mathbb{N}$ and every measurable function 
    $\phi \colon \Vcal \to \mathbb{R}$ satisfying $\|\phi\|_V < \infty$.  
  \end{prop}

  This proposition is a version of Harris' ergodic theorem~\cite{harris1956existence}. 
  For a more contemporary treatment, see~\cite[Section 15]{meyn2012markov}.
  The particular formulation used here follows~\cite{hairer2011yet}.

  To study the regular variation of $\pi$, recall that for a subset 
  $A \subseteq \Vcal$, the gauge function 
  $\rho_A \colon \Vcal \to [0, +\infty]$ is defined by  
  \begin{equation*}  
    \rho_A(\xbf) = \inf\{r > 0: \: r \xbf \in A\}.  
  \end{equation*}  
  Here and throughout the article, the infimum of the empty set is understood 
  to be $+\infty$.  
  
  We say that $K\subseteq \Vcal$ is radially increasing if 
  $sK \supseteq tK$ for $s<t$. Equivalently,
  \[
    \xbf\in K \ \Rightarrow\ s \xbf \in K \quad \text{for all } s\ge 1.
  \]
  We proceed under the following additional assumption.

  \begin{as}\label{as:2:gague}
    \begin{description}  
      \item [$(i)$] \textbf{Monotonicity:} There exists a radially increasing open 
      set $K \subseteq \Vcal$ such that $0 \notin \operatorname{cl}(K)$. Moreover, 
      whenever $s \leq t$, for every 
      $\xbf \in \rho_{K}^{-1}[\{s\}]$ and $\ybf \in \rho_{K}^{-1}[\{t\}]$,  
      \begin{equation*}  
        \Prm(\ybf, rK) \leq \Prm(\xbf, rK)  
      \end{equation*}  
      for every $r > 0$, where
      \begin{equation*}
        \rho_K^{-1}[\{s\}] = \{\xbf\in \Vcal: \: \rho_K(\xbf) = s\}.
      \end{equation*}
      \item [$(ii)$] \textbf{Contractivity:} We have
      \begin{equation*}  
        \rho_* = \sup_{\|\theta\| = 1}  
        \Ebb \left[  
        \left\|\Zbf_1^\theta\right\|^{\alpha} \right] < 1.  
      \end{equation*}  
    \end{description}
  \end{as}

  \begin{figure}
    \begin{center}
      \begin{tikzpicture}[scale=0.7,>=Latex, line cap=round, line join=round]
        \def\angA{10}
        \def\angB{60}
        \def\Rmax{8}

        \def\s{1.6}
        \def\t{2.6}
        \def\r{3.4}   

        \def\angx{30}
        \def\angy{42}

        \path[fill=gray!15]
          (\angA:\r) arc[start angle=\angA,end angle=\angB,radius=\r] --
          (\angB:\Rmax) arc[start angle=\angB,end angle=\angA,radius=\Rmax] -- cycle;

        \draw[thick]
         (\angA:\r) arc[start angle=\angA,end angle=\angB,radius=\r];

        \draw[thick]
          (\angA:\r) -- (\angA:\Rmax);
          \draw[thick]
          (\angB:\r) -- (\angB:\Rmax);

        \draw[thick,dashed]
          (\angB:\Rmax) arc[start angle=\angB,end angle=\angA,radius=\Rmax];

        \node at (35:4.5) {$rK$};

        \draw[thick]
          (\angA:\s) arc[start angle=\angA,end angle=\angB,radius=\s];

        \draw[thick]
          (\angA:\t) arc[start angle=\angA,end angle=\angB,radius=\t];

        \node[left] at (63:\s) {$\rho_K^{-1}(\{t\})$};
        \node[left] at (63:\t) {$\rho_K^{-1}(\{s\})$};

        \coordinate (X) at (\angx:\s);
        \coordinate (Y) at (\angy:\t);

        \fill (X) circle (0.07) node[below left] {$y$};
        \fill (Y) circle (0.07) node[below left] {$x$};
      \end{tikzpicture}
    \end{center}
    \caption{The depiction of sets that appear in the Assumption~\ref{as:2:gague}}
  \end{figure}

  The monotonicity assumption can equivalently be reformulated as follows. 
  If $s < t$, then for any probability measures $\mu$ and $\nu$ 
  concentrated on $\rho_{K}^{-1}[\{s\}]$ and $\rho_{K}^{-1}[\{t\}]$, respectively, 
  \begin{equation}
    \int \Prm(\zbf, rK) \,\nu(\ud \zbf) 
    \leq \int \Prm(\zbf, rK) \,\mu(\ud \zbf)
  \end{equation}
  for every $r > 0$.

  The level set $\rho_K^{-1}[\{s\}]$ can be viewed as a radial boundary of 
  the set $s^{-1}K$. Thus, the monotonicity condition expresses the fact that 
  the farther the starting point of the chain is from the origin, the greater 
  the probability that the chain remains away from the origin after one step.
  
  The second condition in Assumption~\ref{as:2:gague}
  is of a more technical nature and will be used to control $\pi(tK)$.
  To explain its role, recall that $\Krm(\xbf,\cdot) = \Pbb[\Zbf_1^\xbf \in \cdot]$
  denotes the transition kernel of $\{\Zbf_n^\xbf\}_{n \in \mathbb{N}}$.
  The homogeneity of $\Zbf$ implies that
  \begin{equation*}
    \Ebb\left[ \left. \|\Zbf_n^\xbf\|^\alpha \,\right|\, \Zbf_{n-1}^\xbf \right]
      = \|\Zbf_{n-1}^\xbf\|^\alpha \int \|\ybf\|^\alpha
      \,\Krm \left(\frac{\Zbf_{n-1}^\xbf}{\|\Zbf_{n-1}^\xbf\|}, \ud \ybf\right).
  \end{equation*}
  Under the contractivity assumption, this yields
  \begin{equation*}
    \Ebb\left[ \left. \|\Zbf_n^\xbf\|^\alpha \,\right|\, \Zbf_{n-1}^\xbf \right]
    \leq \rho_* \|\Zbf_{n-1}^\xbf\|^\alpha.
  \end{equation*}
  Iterating this inequality gives
  \begin{equation}\label{eq:2:boundZbf}
    \Ebb\left[ \left\|\Zbf_n^\theta\right\|^\alpha \right]
      \leq \left(\rho_*\right)^n,
      \qquad \theta \in \Scal.
  \end{equation}
  Since $0 \notin \operatorname{cl}(K)$, the set $K$ is bounded away from the origin. 
  Hence there exists $\epsilon > 0$, depending on $K$, such that 
  $\rho_{\operatorname{cl}(K)}(\xbf) \geq \epsilon \|\xbf\|^{-1}$ for all $\xbf \in \Vcal$.
  Consequently, for all $\theta \in \Scal$,
  \begin{equation*}
    \Ebb \left[ \rho_{\operatorname{cl}(K)}\!\left(\Zbf_n^\theta\right)^{-\alpha} \right] 
    \leq \left(\rho_*\right)^n / \epsilon^{\alpha}.
  \end{equation*}
  Since $\rho_* < 1$, there exists $k_0 \in \Nbb$ such that
  \begin{equation}\label{eq:2:contr2.0}
    \sup_{\theta \in \Scal} 
    \Ebb \left[ \rho_{\operatorname{cl}(K)}\!\left(\Zbf_{k_0}^\theta\right)^{-\alpha} \right] 
    < 1.
  \end{equation}
  This is the form of the contractivity condition that we will use to control 
  $\pi(tK)$ in the proof of our main result.
        
  We now state and prove a lemma closely following~\cite[Lemma~3]{grey1994regular}.

  \begin{lem}\label{lem:4:supsol}
    Let $k_0$ be a sufficiently large integer such that~\eqref{eq:2:contr2.0} holds.
    There exist a random vector $\Wbf$ and a constant $L>0$ such that, as $t\to\infty$,
    \begin{equation*}
      t \Pbb[\Wbf \in a_tK] \to L,
    \end{equation*}
    and, for every $t>0$,
    \begin{equation*}
      \Pbb[\Wbf \in a_t K] \geq \Ebb\left[\Prm^{k_0}(\Wbf, a_t K)\right],
    \end{equation*}
    where $\Prm^{k_0}$ denotes the $k_0$-step transition kernel of $\Prm$.
  \end{lem}

  \begin{proof}
    To simplify the notation, by considering the $k_0$-step Markov chain 
    $\{\Xbf_{k_0n}^\xbf\}_{n \in \Nbb}$, we can assume without loss of 
    generality that $k_0=1$. Next, consider any $\Vcal$-valued random 
    vector $\Ybf$ that is $\Rvrm(\alpha, \mu_\Ybf)$, such that
    \begin{equation*}
      t \Pbb[\| \Ybf\| > a_t] \to L_\Ybf > 0.
    \end{equation*}
    Then, according to Lemma~\ref{lem:3:key} by approximating $\mathbf{1}_K$ 
	by $\mathcal{C}_0$ functions
	\begin{equation*}
      \limsup_{t \to \infty} t \Ebb\left[\Prm(\Ybf, a_tK)\right] \leq 
      \Ebb \left[\nu_{\Ybf}\left(\mathrm{cl}(K)\right) \right] +
      L_\Ybf \int \Krm\left(\xbf,\mathrm{cl}(K)\right) \mu_\Ybf(\ud \xbf).
	\end{equation*}
	Assumption~\ref{as:2:gague} allows us to control the second term. 
    Indeed, since $\Ybf$ is 
	also $\Rvrm(\alpha,\Theta_\Ybf)$ then
	\begin{equation*}
      \int \Krm\left(\xbf, \mathrm{cl}(K) \right) \mu_\Ybf(\ud \xbf) = 
      \int_\Scal \int_0^{+\infty} 
      \Krm(r \theta,\mathrm{cl}(K)) \alpha r^{-\alpha -1}\ud r 
      \Pbb[\Theta_\Ybf \in \ud \theta].
	\end{equation*}
    Now note that for any $c>0$ and 
	$\xbf \in \Vcal$, $\Krm(c\xbf,\ud \ybf) = \Krm(\xbf,c^{-1} \ud \ybf)$.
	and thus we can continue the last display 
	with
	\begin{equation*}
      \int_\Scal \int_0^{+\infty} 
      \Krm\left(\theta,r^{-1} \mathrm{cl}(K)\right) \alpha r^{-\alpha -1} \ud r 
      \Pbb[\Theta_\Ybf \in \ud \theta].
	\end{equation*}
	The inner integral is bounded by
	\begin{multline*}
      \int_\Vcal \int_0^{+\infty} 1_{\left\{ \xbf \in r^{-1} \mathrm{cl}(K)\right\}} 
      \alpha r^{-\alpha -1} \ud r \Krm(\theta,\ud \xbf)
      = \int_\Vcal \rho_{\operatorname{cl}(K)}(\xbf)^{-\alpha} \Krm(\theta,\ud \xbf) \\ 
      \leq q = \sup_{\theta \in \mathcal{S}} 
      \Ebb\left[\rho_{\operatorname{cl}(K)}(\Zbf_{k_0}^\theta)^{-\alpha}\right].
	\end{multline*}
	This yields
	\begin{equation*}
      \limsup_{t \to \infty} t \Ebb\left[\Prm(\Ybf, a_tK)\right] \leq 
      \Ebb \left[\nu_{\Ybf}\left(\mathrm{cl}(K)\right) \right] + L_\Ybf q.
	\end{equation*}
	Since $q< 1$, we can increase $L_\Ybf$ while keeping the expectation bounded. 
    Therefore, without loss of generality, we can assume that the right-hand 
    side of the last display is strictly smaller than $L_\Ybf$.
	This means that for some $t_0>0$ we have  
	\begin{equation*}
      \Ebb \left[ \Prm(\Ybf, a_tK)  \right] < \Pbb[\Ybf \in a_t K]
	\end{equation*}
	for all $t>t_0$. 
	Define $\Wbf$ to have the conditional distribution of $\mathbf{Y}$ 
    given that $\mathbf{Y} \in a_{t_0}K$.
    We will check that such $\Wbf$ satisfies the claim of the lemma. 
    For $t>t_0$ we have
	\begin{equation*}
      \Ebb\left[\Prm(\Wbf, a_tK) \right] \leq 
      \frac{\Ebb\left[ \Prm\left( \Ybf, a_tK\right) \right]}{\Pbb[ \Ybf \in a_{t_0}K]} 
      < \frac{\Pbb[ \Ybf \in a_{t}K]}{\Pbb[ \Ybf \in a_{t_0}K]} =  \Pbb[\Wbf \in a_tK].
	\end{equation*}
	For $t\leq t_0$ we have, since $K$ is radially increasing, 
	\begin{equation*}
      \Pbb\left[ \Wbf \in a_t K \right] = 
      \Pbb \left[ \left. \Ybf \in a_t K \right| \Ybf \in a_{t_0}K \right] = 1
	\end{equation*}
	so the claimed inequality holds trivially.
  \end{proof}

  We can now state and prove our main result. Recall that $\pi$ denotes the 
  stationary distribution of the process $\{\Xbf_n^\xbf\}_{n\in \Nbb}$. 
  Thus,
  \begin{equation*}
    \pi(\cdot) = \int_{\Vcal} \Prm(\xbf, \cdot)\, \pi(\ud \xbf).
  \end{equation*}

  \begin{thm}\label{thm:2:pi}
    Suppose that Assumptions~\ref{as:2:first},~\ref{as:2:mix} 
    and~\ref{as:2:gague} hold. Assume additionally that, for some constant 
    $C_V > 0$,
    \begin{equation}\label{eq:2:veq}
      \inf_{\|\xbf\|\geq a_t} V(\xbf) \ge C_V t
    \end{equation}
    for all $t\geq 1$, and that $\pi(K) > 0$.
    Then $\pi(\cdot \cap K)$ belongs to $\mathrm{Rv}(\alpha, \mu_\infty)$, 
    where the measure $\mu_\infty$ is determined by
    \begin{equation*}
      \int_K f\left(\ybf\right) \mu_\infty(\ud \ybf)
      = L_\infty\sum_{j=0}^\infty \int_\Vcal\int_\Vcal 
      \Ebb\left[f(\Zbf_j^\ybf) 
      \mathbf{1}_{\left\{ \Zbf_j^\ybf \in K \right\}} \right] 
      \nu_{\xbf}(\ud \ybf) \pi(\ud \xbf)
    \end{equation*}
    for all $f \in \mathcal{C}_0$, where $L_\infty$ is the normalising constant such that $\mu_\infty(\Bcal^c)=1$.
  \end{thm}

\begin{proof}
	We first establish that
	\begin{equation}\label{eq:4:Pfpiclaim1}
		\limsup_{t \to \infty} t \pi (a_tK) <\infty.
	\end{equation}
	To achieve this we will approach $\pi$ with an appropriate choice of the initial condition.
	Namely, we will consider our $k_0$-step Markov chain, with the law of the random 
    vector $\mathbf{W}$ from Lemma~\ref{lem:4:supsol} being the initial distribution.
	More precisely, let $\Xbf^\Wbf = \{ \Xbf_{k_0n}^\Wbf \}_{n \in \Nbb}$ with 
    $\Xbf_0^\Wbf = \Wbf$.
	We next show monotonicity of $\Xbf^\Wbf$ in an appropriate sense. 
	To verify the mentioned monotonicity, note that
	\begin{equation*}
		\left\{ \Wbf \in tK \right\} = \left\{ \rho_{K}(\Wbf)< 1/t   \right\} = 
        \left\{ 1/\rho_{K}(\Wbf)> t   \right\}
	\end{equation*}
	Using the monotonicity condition from Assumption~\ref{as:2:gague}
    	we see that the sequence of random variables $1 / \rho_{K}(\Xbf_n^\Wbf)$ 
	is stochastically decreasing. Indeed, the claim of the 
    aforementioned Lemma~\ref{lem:4:supsol}
	gives
	\begin{align*}
		\Pbb\left[ 1 / \rho_{K}\left(\Xbf_0^\Wbf\right) \geq t \right] & =
		\Pbb\left[ \Xbf_0^\Wbf \in tK \right] \\ & \geq 
		\Pbb\left[ \Xbf_{k_0}^\Wbf \in tK \right] =\Pbb\left[ 1 / 
        \rho_{K}\left(\Xbf_{k_0}^\Wbf\right) \geq t \right].
	\end{align*}
	From here we argue inductively. If for $n \in \Nbb$, $1 / \rho_{K}(\Xbf_{k_0(n-1)}^\Wbf)$ is stochastically
    greater than $1 / \rho_{K}(\Xbf_{k_0n}^\Wbf)$ in the standard stochastic order, 
    then there exists a coupling $(X, Y)$ of 
    $(\rho_{K}(\Xbf_{(n-1)k_0}^\Wbf), \rho_{K}(\Xbf_{nk_0}^\Wbf))$ 
    such that $X \leq Y$ almost surely.  
	Let $\zeta( \rho_{K}(\Xbf_{nk_0}^\Wbf), \cdot)$ be the regular conditional distribution
	of $\Xbf_{nk_0}^\Wbf$ given $\rho_{K}(\Xbf_{nk_0}^\Wbf)$. Then
	\begin{equation*}
		\Ebb \left[\left.\Prm^{k_0}\left(\Xbf_{nk_0}^\Wbf, tK\right) \right| 
        \rho_{K}\left(\Xbf_{nk_0}^\Wbf\right)  \right] 
			= \int \Prm^{k_0}(\zbf, tK) \zeta(\rho_{K}\left(\Xbf_{nk_0}^\Wbf\right), \ud \zbf).
	\end{equation*}
	Next write $\eta$ for the conditional distribution of $\Xbf_{(n-1)k_0}^\Wbf$ given 
	$\rho_{K}(\Xbf_{(n-1)k_0}^\Wbf)$. Then
	\begin{equation*}
		\Ebb \left[\left. \Prm^{k_0}\left(\Xbf_{(n-1)k_0}^\Wbf, tK\right) \right| \rho_{K} 
        \left(\Xbf_{(n-1)k_0}^\Wbf\right)  \right] 
			= \int \Prm^{k_0}(\zbf, tK) \eta(\rho_{K}(\Xbf_{(n-1)k_0}^\Wbf), \ud \zbf).
	\end{equation*}
	Using the monotonicity condition we can write
	\begin{multline*}
		\Pbb \left[ 1 / \rho_{K}(\Xbf_{(n+1)k_0}^\Wbf) \geq t \right] = 
		\Ebb \left[ \Ebb \left[\left. 
				\Prm^{k_0}\left(\Xbf_{nk_0}^\Wbf, tK\right) \right| 
                \rho_{K}\left(\Xbf_{nk_0}^\Wbf\right)  \right] 
		\right] =\\ 
		\Ebb \left[ \int \Prm^{k_0}(\zbf, tK) 
        \zeta(\rho_{K}\left(\Xbf_{nk_0}^\Wbf\right), \ud \zbf) \right]
		= \Ebb\left[ \int \Prm^{k_0}(\zbf, tK) \zeta(Y, \ud \zbf) \right] \\
		\leq \Ebb\left[ \int \Prm^{k_0}(\zbf, tK) \eta(X, \ud \zbf) \right]
		=\Pbb_\xbf[1 / \rho_{K}\left(\Xbf_{nk_0}^\Wbf\right)\geq t].
	\end{multline*}
	Thus the sequence $\left\{\P[\Xbf_{n k_0}^\Wbf \in t K]\right\}_{n \in \mathbb{N}}$ 
	is non-increasing and~\eqref{eq:4:Pfpiclaim1} follows since
	\begin{equation*}
		t\pi(a_t K) \leq \limsup_{n \to \infty} t\Pbb[\Xbf_{nk_0}^\Wbf \in a_tK ] \leq 
		t \Pbb[\Wbf \in a_tK].
	\end{equation*}
 
	From here we go back to the Markov chain started at some fixed $\xbf \in \Vcal$.
	For $f \in \Ccal_0$ denote $f_K(x) = f(x)\mathbf{1}_{\{ x \in K \}}$. 
    For $n \in \Nbb$ write
	\begin{multline*}
		t \int f_K\left(a_t^{-1} \ybf \right) \pi(\ud \ybf) = \\
		t \left(\int f_K\left(a_t^{-1} \ybf \right) \pi(\ud \ybf) - 
		 \Ebb\left[ f_K \left( a_t^{-1} \Xbf_n^\xbf  \right) \right] \right)
		+t \Ebb \left[ f_K \left( a_t^{-1} \Xbf_n^\xbf  \right) \right].
	\end{multline*}
	The first term is bounded by
	\begin{equation*}
		C_{\pi} t \rho^n (1+V(\xbf)) \left\| f_K(a_t^{-1} \cdot) - 
		\int f_K\left(a_t^{-1} \ybf \right) \pi(\ud \ybf)
		\right\|_V
	\end{equation*}
	one only needs to argue that the term involving $\|\cdot\|_V$ is $O(1/t)$.
	The integral is easily manageable using~\eqref{eq:4:Pfpiclaim1}, 
    while $\|f_K(a_t^{-1} \cdot)\|_V$ 
	can be controlled using~\eqref{eq:2:veq}.
	Now analyse the second term by first taking $t \to \infty$ followed by $n \to \infty$. 
	We have
	\begin{equation}\label{eq:2:16}
		\lim_{t \to \infty}
		t \Ebb \left[ f_K \left( a_t^{-1} \Xbf_n^\xbf  \right) \right]= 
		\sum_{j=0}^{n-1} \Ebb \left[ \int_\Vcal 
        \Ebb[f_K(\Zbf_j^\ybf) ] \nu_{\Xbf_{n-j-1}^\xbf } (\ud y) \right].
	\end{equation}
	We now argue that the sum converges as $n \to \infty$. 
    	Recall that there exists $\epsilon$ such that $\epsilon \Bcal \cap K =\emptyset$. Then
    	for non-negative integers $n$ and $j <n$ we have
	\begin{align*}
		\Ebb \left[ \int_\Vcal \Ebb[f_K(\Zbf_j^\ybf) ] 
        \nu_{\Xbf_{n-j-1}^\xbf } (\ud y) \right]& \leq 
		\Ebb \left[ \int_\Vcal c \Pbb[ \|\Zbf_j^\ybf\| \geq \epsilon ] 
        \nu_{\Xbf_{n-j-1}^\xbf } (\ud y) \right] \\
		& \leq c \epsilon^{-\alpha} \left( \rho_* \right)^j 
        \mathbb{E}\left[\|\Xbf_{n-j-1}^\xbf\|^\alpha\right]
	\end{align*}
    for some $c>0$.
	This allows us to control the terms in the sum~\eqref{eq:2:16} for big values of $j$. 
	If $j$ is fixed then one just needs to note that
	\begin{equation*}
		\xbf \mapsto F_j(\xbf) = \int_\Vcal \Ebb[f_K (\Zbf_j^\ybf)] \: \nu_\xbf(\ud \ybf) 
	\end{equation*}
	satisfies $\|F_j \|_V <\infty$. Thus for fixed $j$,
	\begin{equation*}
		\Ebb \left[ \int_\Vcal \Ebb[f_K(\Zbf_j^\ybf) ] \nu_{\Xbf_{n-j-1}^\xbf } (\ud y) \right] 
		= \Ebb \left[ F_j \left( \Xbf_{n-j-1}^\xbf \right) \right] \to \int_\Vcal F_j(\ybf) \pi(\ud \ybf).
	\end{equation*}
	The last two remarks allow one to conclude that
	\begin{equation*}
		\lim_{n \to \infty}
		\sum_{j=0}^{n-1} \Ebb \left[ \int_\Vcal \Ebb[f_K(\Zbf_j^\ybf) ] \nu_{\Xbf_{n-j-1}^\xbf } (\ud y) \right]
		=\sum_{j=0}^\infty \int_\Vcal\int_\Vcal \Ebb[f_K(\Zbf_j^\ybf) ] \nu_{\xbf } (\ud y) \pi(\ud \xbf).
	\end{equation*}
	This yields 
	\begin{equation*}
		\lim_{t \to \infty}t \int f_K\left(a_t^{-1} \ybf \right) \pi(\ud \ybf)
		=\sum_{j=0}^\infty \int_\Vcal\int_\Vcal \Ebb[f_K(\Zbf_j^\ybf) ] \nu_{\xbf } (\ud y) \pi(\ud \xbf)
	\end{equation*}
	This establishes multivariate regular variation of $\pi( \cdot \cap K)$. 
\end{proof}

\section{Random coefficient linear models}\label{sec:ar}

  We now turn to random coefficient linear models. 
  We first introduce the general framework and then present our results 
  in two distinct settings.

\subsection{A general setting}

  Random coefficient linear models are typically formulated in 
  finite-di\-men\-sion\-al spaces. 
  However, the arguments developed in this work allow us to work in a more general 
  setting in which the underlying space may be infinite-dimensional. Since the random 
  series considered below involve products, we extend the previous framework by assuming 
  that the Banach space $\Vcal$ is equipped with a continuous, bilinear, and associative 
  multiplication, making it a Banach algebra (for a comprehensive introduction to Banach 
  algebras, see, e.g., \cite{conway2019course}). Thus, throughout this section, $\Vcal$ 
  denotes a Banach algebra.

  We study random series of the form
  \begin{equation*}
    \Sbf = \sum_{j=1}^\infty \Mbf_j \Nbf_j 
  \end{equation*}
  where $\Mbf_j$ and $\Nbf_j$ are random elements of $\Vcal$.
  We assume throughout that the sequence $(\Nbf_j)_{j \ge 1}$ is i.i.d. and call it innovations.
  We denote a generic copy by $\Nbf$ and assume that $\Nbf$ is regularly varying.

  Our objective is to determine how the dependence structure between the sequences 
  $(\Mbf_j)_{j \ge 1}$ and $(\Nbf_j)_{j \ge 1}$ affects the tail asymptotics of 
  the series. To this end, we consider two distinct dependence regimes.

  \begin{description}
    \item[Predictable case:] There exists a filtration $(\mathcal{F}_j)_{j \ge 0}$ such that 
      $\Mbf_{j+1}$ and $\Nbf_j$ are $\mathcal{F}_j$-measurable, and 
      $\mathcal{F}_j$ is independent of $\sigma(\Nbf_{j+1}, \Nbf_{j+2}, \ldots)$.

    \item[Adapted case:] There exists a filtration $(\mathcal{F}_j)_{j \ge 0}$ such that 
      $\Mbf_j$ and $\Nbf_j$ are $\mathcal{F}_j$-measurable for every $j \ge 1$, and 
      $\mathcal{F}_j$ is independent of $\sigma(\Nbf_{j+1}, \Nbf_{j+2}, \ldots)$.
  \end{description}

  To describe how large innovations propagate through the coefficient process, we impose 
  additional structural assumptions. Specifically, in both cases, we assume that 
  $(\Mbf_j)_{j \ge 1}$ is a Markov chain generated by iterating a deterministic function 
  $\psi \colon \Vcal \times \Vcal \to \Vcal$ satisfying the condition stated below. 
  In particular, the coefficient chains in the two cases have the same transition kernel 
  and, when initialized with the same distribution, the same law.
	
  \begin{as}\label{as:psi}
    There exists a constant $C_\psi > 0$ such that for all $\mbf, \nbf \in \Vcal$,
    \begin{equation}\label{sass1}
      \|\psi(\mbf, \nbf)\| \le C_{\psi} \|\mbf\| \bigl( \|\nbf\| + 1 \bigr).
    \end{equation}
    Moreover, there exist functions 
    $\chi, \zeta \colon \Vcal \times \Vcal \to \Vcal$ 
    such that uniformly on bounded sets 
    \begin{equation}\label{sass2}
      \lim_{t \to \infty} t^{-1} \psi(t \mbf, \nbf) = \chi(\mbf, \nbf),
      \qquad
      \lim_{t \to \infty} t^{-1} \psi(\mbf, t \nbf) = \zeta(\mbf, \nbf),
    \end{equation}
    for $\mbf, \nbf \in \Vcal$,
    where $\chi$ is continuous in its first argument and $\zeta$ is continuous 
    in its second argument.
  \end{as}

  We impose the following condition on the distribution of the innovations.

  \begin{as}\label{as:N}
    The law of $\Nbf$ is regularly varying with index $\alpha$, 
    denoted by $\mathrm{RV}(\alpha, \Theta_{\Nbf})$, and
    \begin{equation*}
      \lim_{M \to \infty} \, \limsup_{t \to \infty} 
      \frac{\Pbb\bigl( \|\Nbf_1\| \, \|\Nbf_2\| > t,\; 
      \|\Nbf_1\| \in (M, t/M] \bigr)}
      {\Pbb\bigl( \|\Nbf\| > t \bigr)} 
      = 0.
    \end{equation*}
  \end{as}

\subsection{The adapted case}

  We now specify the framework for the adapted case. 
  Let $(\Nbf_j)_{j \ge 1}$ be a sequence of i.i.d.\ $\Vcal$-valued random 
  elements. Starting from a deterministic initial condition 
  $\Mbf_0 = \mbf \in \Vcal$, define $(\Mbf_j)_{j \ge 1}$ recursively by
  \begin{equation*}
    \Mbf_j = \psi(\Mbf_{j-1}, \Nbf_j), 
    \qquad j \ge 1,
  \end{equation*}
  where $\psi \colon \Vcal \times \Vcal \to \Vcal$ satisfies 
  Assumption~\ref{as:psi}. We study the random series
  \begin{equation*}
    \Sbf = \sum_{j=1}^\infty \Mbf_j \Nbf_j.
  \end{equation*}

\begin{example}
Let $\Vcal=\mathbb{R}$ and let $(\Nbf_j)_{j\geq1}$ be an i.i.d.
sequence of heavy-tailed innovations. Fix $c>0$ and define
\begin{equation*}
\Mbf_0=1,
\qquad
\Mbf_j=c\Mbf_{j-1}|\Nbf_j|,
\qquad j\geq1.
\end{equation*}
Then
\begin{equation*}
\Sbf
=
\sum_{j=1}^{\infty}\Mbf_j\Nbf_j
=
\sum_{j=1}^{\infty}
c^j\left(\prod_{k=1}^j|\Nbf_k|\right)\Nbf_j.
\end{equation*}
Setting $A_j=c|\Nbf_j|$ and $B_j=A_j\Nbf_j$, we see that, whenever
the series converges,
\begin{equation*}
\Sbf
\stackrel{d}{=}
A\widehat{\Sbf}+B
=
c|\Nbf|\bigl(\widehat{\Sbf}+\Nbf\bigr),
\end{equation*}
where $\widehat{\Sbf}$ is an independent copy of $\Sbf$. Thus $\Sbf$
is the stationary solution of an affine stochastic recursion, belonging
to the classical class of perpetuities; see, for example,
\cite{goldie1991implicit}.

This model is adapted because $\Mbf_j$ depends on the current innovation
$\Nbf_j$. In particular, the same shock determines both the random
coefficient $A_j$ and the additive term
$B_j=c|\Nbf_j|\Nbf_j$. A large innovation is therefore amplified
quadratically, illustrating the mechanism behind the change of the tail
index from $\alpha$ to $\alpha/2$ in the adapted case.
\end{example}

  To analyze the tail behavior of $\Sbf$, we first consider its partial sums. 
  For $n \in \mathbb{N}$, let
  \begin{equation*}
    \Sbf_n = \sum_{j=1}^n \Mbf_j \Nbf_j.
  \end{equation*}

  Before stating the formal results, we provide the main heuristic underlying 
  our approach. Observe that
  \[
    \Mbf_j \Nbf_j = \psi(\Mbf_{j-1}, \Nbf_j)\, \Nbf_j,
  \]
  where $\psi$ is asymptotically homogeneous in its second argument. Hence, 
  when $\|\Nbf_j\|$ is large,
  \begin{equation*}
    \Mbf_j \Nbf_j
    \approx \|\Nbf_j\|^2
    \zeta\left(\Mbf_{j-1},\frac{\Nbf_j}{\|\Nbf_j\|}\right)
    \frac{\Nbf_j}{\|\Nbf_j\|}.
  \end{equation*}
  This suggests that $\Mbf_j\Nbf_j$ is regularly varying with index 
  $\alpha/2$, with the limiting angular contribution determined by
  \begin{equation*}
    \Ubf
    = \zeta(\Mbf_{j-1},\Theta_\Nbf)\Theta_\Nbf.
  \end{equation*}
  A large value of $\|\Nbf_j\|$ also propagates to the next coefficient. Indeed,
  \begin{equation*}
    \Mbf_{j+1}\Nbf_{j+1}
    \approx \|\Nbf_j\|
    \chi\left(\zeta(\Mbf_{j-1},\Theta_\Nbf),\Nbf_{j+1}\right)
    \Nbf_{j+1}
    = o_{\Pbb}\left(\|\Nbf_j\|^2\right).
  \end{equation*}
  Thus, at the scale $\|\Nbf_j\|^2$, all subsequent terms are negligible. More 
  precisely, the preceding discussion suggests the conditional convergence
  \begin{equation*}
    \mathcal{L}\left(
      \left.
      \left\{
        \|\Nbf_j\|^{-2}\Mbf_{j+k}\Nbf_{j+k}
      \right\}_{k\in\mathbb{N}_0}
      \,\right|\, \|\Nbf_j\|>t
    \right)
    \Longrightarrow
    \mathcal{L}\left(\Ubf,0,0,\ldots\right)
  \end{equation*}
  as $t\to\infty$, in the sense of finite-dimensional distributions.
  Consequently, on the event that $\Nbf_j$ is the single large innovation, 
  the partial sum is asymptotically dominated by its $j$th summand:
  \begin{equation*}
    \Sbf_n \approx \|\Nbf_j\|^2\Ubf.
  \end{equation*}
  The lemma proposition below make this heuristic precise.

  \begin{lem}\label{lem:one-large-innovation}
    Suppose that \eqref{sass1} holds and that $\Nbf$ is regularly varying 
    with index $\alpha>0$. Let $(a_t)$ satisfy
    \begin{equation*}
      t\Pbb[\|\Nbf\|>a_t]\longrightarrow 1.
    \end{equation*}
    Then, for every $n\in\mathbb{N}$ and every $f\in\Ccal_0^u$,
    \begin{equation}\label{eq:one-large-innovation}
      t\left|
        \Ebb\left[f\left(a_t^{-2}\Sbf_n\right)\right]
        -
        \sum_{k=1}^n
        \Ebb\left[
          f\left(a_t^{-2}\Mbf_k\Nbf_k\right)
        \right]
      \right|
      \longrightarrow 0.
    \end{equation}
  \end{lem}

  \begin{proof}
    Since $\Nbf$ is regularly varying with index $\alpha$, the sequence
    $(a_t)$ may be chosen to be regularly varying with index $1/\alpha$. Write
    \begin{equation*}
      \Tbf_k=\Mbf_k\Nbf_k,
      \qquad k\geq1,
    \end{equation*}
    and set $R_k=1+\|\Nbf_k\|$.
    Iterating \eqref{sass1}, we obtain
    \begin{equation*}
      \|\Mbf_k\|
      \leq C_\psi^k\|\Mbf_0\|
      \prod_{\ell=1}^kR_\ell.
    \end{equation*}
    Consequently, for some deterministic constant $c_k>0$,
    \begin{equation}\label{eq:Tk-product-bound}
      \|\Tbf_k\|
      \leq R_k^2 Q_k, \quad 
      Q_k=c_k\prod_{\ell=1}^{k-1}R_\ell. 
    \end{equation}
    Put $\beta=\alpha/2$ and choose $\delta \in (0,\alpha/2)$.
    Then $Q_k$ is independent of $R_k$ and satisfies
    \begin{equation*}
      \Ebb[Q_k^{\beta+\delta}]<\infty.
    \end{equation*}
    Moreover, $R_k^2$ is regularly varying with index $\beta$. Hence,
    by Breiman's lemma and \eqref{eq:Tk-product-bound}, for every $u>0$,
    \begin{equation*}\label{eq:Tk-tail-bound}
      \limsup_{t\to\infty}
      t\Pbb[\|\Tbf_k\|>u a_t^2]
      \leq u^{-\beta}\Ebb[Q_k^\beta]
      <\infty.
    \end{equation*}
    In particular,
    \begin{equation*}\label{eq:Tk-tail-tightness}
      \lim_{M\to\infty}\limsup_{t\to\infty}
      t\Pbb[\|\Tbf_k\|>M a_t^2]=0.
    \end{equation*}

    We next show that two distinct $\Tbf_k$'s cannot simultaneously contribute
    at the scale $a_t^2$. Fix $1\leq i<j\leq n$ and choose
    $p\in(\alpha/4,\alpha/3)$.    
      It follows from \eqref{eq:Tk-product-bound} that
    \begin{equation*}
      \|\Tbf_i\|^p\|\Tbf_j\|^p
      \leq c_{i,j}
      \left(\prod_{\ell=1}^{i-1}R_\ell^{2p}\right)
      R_i^{3p}
      \left(\prod_{\ell=i+1}^{j-1}R_\ell^p\right)
      R_j^{2p}
    \end{equation*}
    for some deterministic constant $c_{i,j}>0$. 
    The random variables on the right hand side
    are independent and have finite expectation, 
    since $3p<\alpha$. Therefore, for every
    $u,v>0$,
    \begin{equation*}
      t\Pbb\left[
        \|\Tbf_i\|>u a_t^2,\,
        \|\Tbf_j\|>v a_t^2
      \right] 
      \leq
      \frac{
        t\Ebb[\|\Tbf_i\|^p\|\Tbf_j\|^p]
      }{
        u^pv^pa_t^{4p}
      }
      \longrightarrow0,
    \end{equation*}
    because $(a_t)$ is regularly varying with index $1/\alpha$ and
    $4p>\alpha$. Thus,
    \begin{equation*}\label{eq:pairwise-negligibility}
      t\Pbb\left[
        \|\Tbf_i\|>u a_t^2,\,
        \|\Tbf_j\|>v a_t^2
      \right]
      \longrightarrow0.
    \end{equation*}

    We now pass from pairwise negligibility to the one-large-innovation
    approximation. 
    Since $0\notin\operatorname{supp}(f)$, there exists $r_0>0$ such that $f(\xbf)=0$ whenever $\|\xbf\|\leq r_0$.
    Let $c=r_0/n$.
    Fix $M>c$ and $\eta\in(0,c)$, and define
    \begin{equation*}
      H_t(M)
      =
      \left\{
        \max_{1\leq k\leq n}\|a_t^{-2}\Tbf_{k}\|\leq M
      \right\}.
    \end{equation*}
    On $H_t(M)$, we have
    \begin{equation*}
      \|a_t^{-2}S_n\|\leq nM.
    \end{equation*}

    Since the restriction of $f$ to $nM\Bcal$ is uniformly continuous,
    its modulus of continuity on this set,
    \begin{equation*}
      \omega_{f,nM}(\delta)
      =
      \sup\left\{
        |f(\xbf)-f(\ybf)|:
        \xbf,\ybf\in nM\Bcal,\ 
        \|\xbf-\ybf\|\leq\delta
      \right\},
    \end{equation*}
    satisfies $\omega_{f,nM}(\delta)\to 0$ as $\delta\to 0$.
    For $1\leq k\leq n$, let
    \begin{equation*}
      G_{k,t}(\eta)
      =
      \left\{
        \|a_t^{-2}\Tbf_{k}\|>c,\,
        \max_{j\neq k}\|a_t^{-2}\Tbf_{j}\|\leq\eta
      \right\}.
    \end{equation*}
    On $G_{k,t}(\eta)\cap H_t(M)$
    \begin{equation*}
      \left|
        f(a_t^{-2}S_n)-\sum_{j=1}^nf(a_t^{-2}\Tbf_{j})
      \right| 
      = \left|
        f(a_t^{-2}S_n)-f(a_t^{-2}\Tbf_{k})
      \right| 
      \leq\omega_{f,nM}((n-1)\eta).
    \end{equation*}
    On
    \begin{equation*}
      H_t(M)\setminus\bigcup_{k=1}^nG_{k,t}(\eta),
    \end{equation*}
    the difference
    \begin{equation*}
       f(a_t^{-2}S_n)-\sum_{j=1}^nf(a_t^{-2}\Tbf_{j})   \end{equation*}
    can be nonzero only if there exist distinct indices $i,j$ such that
    \begin{equation*}
      \|a_t^{-2}\Tbf_{i}\|>c
      \qquad\text{and}\qquad
      \|a_t^{-2}\Tbf_{j}\|>\eta.
    \end{equation*}

    Consequently,
    \begin{align*}
      &t\left|
        \Ebb[f(a_t^{-2}S_n)]
        -
        \sum_{k=1}^n\Ebb[f(a_t^{-2}\Tbf_{k})]
      \right| \\
      &\quad\leq
      \omega_{f,nM}((n-1)\eta)
      \sum_{k=1}^n
      t\Pbb[\|\Tbf_k\|>c a_t^2] \\
      &\qquad+
      (n+1)\|f\|_\infty
      \sum_{i\neq j}
      t\Pbb\left[
        \|\Tbf_i\|>c a_t^2,\,
        \|\Tbf_j\|>\eta a_t^2
      \right] \\
      &\qquad+
      (n+1)\|f\|_\infty
      \sum_{k=1}^n
      t\Pbb[\|\Tbf_k\|>M a_t^2].
    \end{align*}
    Taking first the limit superior as $t\to\infty$, then letting
    $\eta\downarrow0$, and finally letting $M\to\infty$, proves \eqref{eq:one-large-innovation}.
  \end{proof}

  \begin{prop}
    Suppose that Assumptions~\ref{as:psi} and~\ref{as:N} hold. Assume additionally that
    \begin{equation*}
      C_n
      :=
      \sum_{k=1}^n
      \Ebb\left[
        \left\|
          \zeta(\Mbf_{k-1},\Theta_\Nbf)\Theta_\Nbf
        \right\|^{\alpha/2}
      \right]
      >0.
    \end{equation*}
    Then, for every $n \in \mathbb{N}$, the law of $\Sbf_n$ belongs to
    $\mathrm{Rv}(\alpha/2,\mu_n)$, where the measure $\mu_n$ is determined by
    \begin{align*}
      \int_{\Vcal_0} f(\xbf)\,\mu_n(\ud\xbf)
      &=
      L_n\sum_{k=1}^n
      \int_0^\infty
      \Ebb\left[
        f\left(
          r^2\zeta(\Mbf_{k-1},\Theta_\Nbf)\Theta_\Nbf
        \right)
      \right]
      \alpha r^{-\alpha-1}\ud r
    \end{align*}
    for all $f\in\mathcal{C}_0$, with
    $L_n=C_n^{-1}$
    In particular, $L_n$ is chosen so that $\mu_n(\Bcal^c)=1$.
  \end{prop}

  \begin{proof}
    Let $(a_t)$ be a scaling sequence such that
    \begin{equation*}
      \lim_{t\to\infty}t\Pbb\bigl[\|\Nbf\|>a_t\bigr]=1.
    \end{equation*}
    We first identify the limiting measure corresponding to the scaling
    sequence $(a_t^2)$. For $f\in\Ccal_0^u$,
    by the merit of~\eqref{eq:one-large-innovation} we only need to analyse 
    \begin{equation*}
        \sum_{k=1}^nf(a_t^{-2}\Mbf_k\Nbf_k).
    \end{equation*}
    For each fixed $k$, the random element $\Mbf_{k-1}$ is independent of
    $\Nbf_k$. Therefore, the asymptotic homogeneity of $\psi$ in its second
    argument and the regular variation of $\Nbf_k$ imply that
    \begin{multline*}
      \lim_{t\to\infty}
      t\Ebb\left[
        f\left(
          a_t^{-2}\psi(\Mbf_{k-1},\Nbf_k)\Nbf_k
        \right)
      \right] \\
      =
      \int_0^\infty
      \Ebb\left[
        f\left(
          r^2
          \zeta(\Mbf_{k-1},\Theta_\Nbf)
          \Theta_\Nbf
        \right)
      \right]
      \alpha r^{-\alpha-1}\ud r,
    \end{multline*}
    where $\Theta_\Nbf$ is independent of $\Mbf_{k-1}$.

    Combining this convergence with
    \eqref{eq:one-large-innovation}, we obtain
    \begin{align*}
      \lim_{t\to\infty}
      t\Ebb\left[f\left(a_t^{-2}\Sbf_n\right)\right]
      &=
      \sum_{k=1}^n
      \int_0^\infty
      \Ebb\left[
        f\left(
          r^2
          \zeta(\Mbf_{k-1},\Theta_\Nbf)
          \Theta_\Nbf
        \right)
      \right]
      \alpha r^{-\alpha-1}\ud r.
    \end{align*}
    The claim follows.
  \end{proof}

  Before proving Theorem~\ref{thm:2:RCLM1}, we record a simple lemma ensuring 
  that the random series $\Sbf$ is well defined. Define the function 
  $\psi^*\colon\Vcal\to[0,+\infty)$ by
  \begin{equation*}
    \psi^*(\nbf)
    =
    \sup_{\mbf\neq0}
    \frac{\|\psi(\mbf,\nbf)\|}{\|\mbf\|}.
  \end{equation*}
  We assume that $\psi^*$ is measurable. Assumption~\ref{as:psi} implies that
  \begin{equation*}
    \psi^*(\nbf)
    \leq C_\psi\bigl(\|\nbf\|+1\bigr),
    \qquad \nbf\in\Vcal.
  \end{equation*}
  In what follows, we impose the contractivity condition
  \begin{equation}\label{eq:2:contr}
    \Ebb\left[\psi^*(\Nbf_1)^{\alpha/2}\right]<1.
  \end{equation}

  \begin{lem}
    Suppose that~\eqref{eq:2:contr} holds. If 
    $\mathbb{E}[\log^+\|\Nbf\|]<\infty$, then the series defining $\Sbf$ 
    converges in $\Vcal$ almost surely.
  \end{lem}

  \begin{proof}
    Since $\Vcal$ is complete, it is sufficient to verify that
    \begin{equation*}
      \sum_{j=1}^\infty \|\Mbf_j \Nbf_j\| <\infty
    \end{equation*}
    almost surely. We have
    \begin{equation}\label{eq:4:estM}
      \|\Mbf_j\| = \|\psi(\Mbf_{j-1}, \Nbf_j)\| \leq 
      \|\Mbf_{j-1}\| \psi^*(\Nbf_j) 
      \leq \|\Mbf_0\| \prod_{k=1}^{j} \psi^*(\Nbf_k).
    \end{equation}
    Under~\eqref{eq:2:contr}, the product on the right-hand side converges 
    to $0$ exponentially fast. The integrability assumption ensures that
    \begin{equation*}
      \limsup_{n \to \infty} \frac{1}{n} \log \|\Nbf_n\| \leq 0 
      \quad \text{a.s.}
    \end{equation*}
    A simple ratio test then shows that
    \begin{equation*}
      \sum_{j=1}^\infty \|\Nbf_j\| 
      \prod_{k=1}^{j} \psi^*(\Nbf_k)<\infty
      \quad \text{a.s.},
    \end{equation*}
    which proves the claim.
  \end{proof}

\begin{thm}\label{thm:2:RCLM1}
	Let the Assumptions~\ref{as:psi} and~\ref{as:N} be in force. Assume additionally that
    \begin{equation*}\label{eq:nondegenerate-limit}
      C_\infty
      :=
      \sum_{k=1}^\infty
      \Ebb\left[
        \left\|
          \zeta(\Mbf_{k-1},\Theta_\Nbf)\Theta_\Nbf
        \right\|^{\alpha/2}
      \right]
      >0.
    \end{equation*}
	If~\eqref{eq:2:contr} holds then
	the law of $\Sbf$ is $\Rvrm(\alpha/2, \mu_\infty)$ with
	\begin{equation*}
    		\int f\left( x \right)\mu_\infty(\mathrm{d}x)
		= L_\infty
        \sum_{k=1}^\infty 
		\int_0^\infty\Ebb\left[f\left(r^2  \zeta(\Mbf_{k-1}, \Theta_\Nbf)\Theta_\Nbf\right)\right] \alpha r^{-\alpha -1}\mathrm{d}r,
	\end{equation*}
	where $L_\infty>0$ is a normalising constant such that $\mu_\infty(\Bcal^c)=1$.
\end{thm}

\begin{proof}
	The argument follows the simple idea of approximating $\Sbf$ via $\Sbf_n$ for large $n \in \mathbb{N}$.
	To achieve that it is sufficient to control the tails of
	\begin{equation*}
		\Sbf_{n, \infty} = \Sbf-\Sbf_n=\sum_{j=n+1}^\infty \Mbf_j\Nbf_j.
	\end{equation*}
	If we show that for any $\varepsilon>0$,
        \begin{equation}\label{eq:2:claim12}
			\lim_{n \to \infty} \limsup_{t \to \infty} t \Pbb\left[\| \Sbf_{n,\infty} \| > \eps a_t^2\right] =0	
	\end{equation}
	then using the uniform continuity of $f$ and the fact that 
        \begin{equation*}
            \lim_{M \to \infty} \lim_{n\to \infty} \limsup_{t\to \infty} t \P \left[\|\Sbf_n\| \ge M a_t^2 \right] = 0
        \end{equation*}
        the main result will follow.
        
        Write
	\begin{equation*}
		\|\Sbf_{n, \infty}\| \leq 
		C_\psi\|\Mbf_n\|\sum_{j=n+1}^\infty (\|\Nbf_j\|+1)^2 \prod_{k=n+1}^{j-1} \psi^*(\Nbf_k).
	\end{equation*}
	By the main result of~\cite{grey1994regular} applied to 
        \begin{equation*}
            \Sbf^*_n = \sum_{j=n+1}^\infty (\|\Nbf_j\|+1)^2 \prod_{k=n+1}^{j-1} \psi^*(\Nbf_k)
        \end{equation*}
        we have
        \begin{equation*}
		\lim_{t\to \infty} t\P [\Sbf^*_n > \eps a_t^2] = \frac{\eps^{-\alpha/2}}{1 - \E [\psi^*(\Nbf)^{\alpha/2}]}
        \end{equation*}
        and
		\begin{equation*}
			\limsup_{t \to \infty} t \Pbb\left[\| \Sbf_{n,\infty} \| > \eps a_t^2\right] 
			\le \frac{C_\psi^{\alpha/2}\eps^{-\alpha / 2} \Ebb[ \| \Mbf_n\|^{\alpha/2}] }{ 1-\Ebb[\psi^*(\Nbf)^{\alpha/2}]}.
		\end{equation*}
		By the estimate~\eqref{eq:4:estM} and~\eqref{eq:2:contr} we have
		\begin{equation*}
			\lim_{n \to \infty} \Ebb\left[ \|\Mbf_n\|^{\alpha/2} \right] =0.
		\end{equation*}
		This secures~\eqref{eq:2:claim12} and proves the result.
	\end{proof}

  \subsection{The predictable case}\label{subsec:2:pred}

  In the predictable case, we consider a random series of the form
  \begin{equation*}
    \Sbf=\sum_{j=1}^\infty\Mbf_j\Nbf_j,
  \end{equation*}
  where $(\Nbf_j)_{j\geq1}$ is a sequence of i.i.d.\ $\Vcal$-valued random
  elements. Starting from a deterministic initial condition
  $\Mbf_1=\mbf\in\Vcal$, define $(\Mbf_j)_{j\geq1}$ recursively by
  \begin{equation*}
    \Mbf_{j+1}=\psi(\Mbf_j,\Nbf_j),
    \qquad j\geq1,
  \end{equation*}
  where $\psi\colon\Vcal\times\Vcal\to\Vcal$ satisfies
  Assumption~\ref{as:psi}. In particular, $\Mbf_j$ is independent of
  $\Nbf_j$ for every $j\geq1$.

  \begin{example}
    Let $d\geq2$, let $\Vcal=M_d(\mathbb{R})$ be equipped with the
    operator norm and the usual matrix multiplication, and denote by
    $\mathbf{e}_1$ the first vector of the canonical basis of $\mathbb{R}^d$.
    Let $(\Zbf_j)_{j\geq1}$ be an iid sequence of $\mathbb{R}^d$-valued
    random vectors.
    We embed the innovations into $\Vcal$ by setting $\Nbf_j=\Zbf_j\mathbf{e}_1^\top$.
    Fix matrices $A,B\in M_d(\mathbb{R})$ and a deterministic positive
    semidefinite matrix $H_1$. Put $\Mbf_1=H_1^{1/2}$
    and consider the intercept-free BEKK-type recursion
    \begin{equation}\label{eq:bekk-recursion}
    \Rbf_j
    =
    H_j^{1/2}\Zbf_j,
    \qquad
    H_{j+1}
    =
    A\Rbf_j\Rbf_j^\top A^\top
    +
    BH_jB^\top.
  \end{equation}
  Here $\Rbf_j\in\mathbb{R}^d$ represents the vector of returns and
  $H_j$ is its conditional covariance matrix.

  To express~\eqref{eq:bekk-recursion} in the framework of this
  section, define, for $m,n\in M_d(\mathbb{R})$,
  \begin{equation*}
    \psi(m,n)
    =
    \left(
      Amn\mathbf{e}_1\mathbf{e}_1^\top n^\top m^\top A^\top
      +
      Bmm^\top B^\top
    \right)^{1/2},
  \end{equation*}
  where the square root denotes the unique positive semidefinite
  square root. BEKK models, introduced by Engle and Kroner~\cite{engle1995multivariate}, constitute a prominent class of multivariate GARCH specifications for describing the joint evolution of conditional variances and covariances of financial returns.
  We have
  \begin{equation*}
    \sum_{j=1}^{\infty}\Mbf_j\Nbf_j = \left(\sum_{j=1}^{\infty}R_j\right)e_1^\top.
  \end{equation*}
   Thus, up to the rank-one embedding $x\mapsto xe_1^\top$, the variable $\Sbf$ represents the cumulative return generated by the volatility process. Since the intercept is zero, 
   the origin is an absorbing state and, 
   under the contractivity assumption, the conditional volatility decays to zero. 
   Consequently, $\Sbf$ may be viewed as the total response of the system, initiated from $H_1$, accumulated before the volatility cascade dies out. 
   The limiting tail measure therefore describes the magnitude and direction of an exceptionally large 
   cumulative response. This interpretation is specific to the contractive zero-intercept model; for a stationary BEKK process with a positive intercept, the unnormalised infinite sum of returns would not generally converge.
  \end{example}

  As in the adapted case, we first study the partial sums
  \begin{equation*}
    \Sbf_n=\sum_{j=1}^n\Mbf_j\Nbf_j,
    \qquad n\in\mathbb{N}.
  \end{equation*}

  Before stating the formal results, we describe the main heuristic behind
  the predictable case. Fix $k\leq n$ and write
  \begin{equation*}
    \widehat{\Theta}_k
    =
    \frac{\Nbf_k}{\|\Nbf_k\|}
  \end{equation*}
  whenever $\Nbf_k\neq0$. Observe that
  \[
    \Mbf_{k+1}\Nbf_{k+1}
    =
    \psi(\Mbf_k,\Nbf_k)\Nbf_{k+1}.
  \]
  When $\|\Nbf_k\|$ is large, the asymptotic homogeneity of $\psi$ in its
  second argument gives
  \begin{equation*}
    \Mbf_{k+1}
    \approx
    \|\Nbf_k\|
    \zeta(\Mbf_k,\widehat{\Theta}_k).
  \end{equation*}
  The asymptotic homogeneity of $\psi$ in its first argument then yields
  \begin{equation*}
    \Mbf_{k+2}
    \approx
    \|\Nbf_k\|
    \chi\left(
      \zeta(\Mbf_k,\widehat{\Theta}_k),
      \Nbf_{k+1}
    \right).
  \end{equation*}
  Thus, a large innovation $\Nbf_k$ propagates through the subsequent
  coefficients and generates a cluster of large summands.

  To describe this cluster, let $\Theta_\Nbf$ have the spectral distribution
  associated with the regular variation of $\Nbf$, independently of
  $\Mbf_k$ and $(\Nbf_{k+j})_{j\geq1}$. Define
  \begin{equation*}
    \Ubf_0^k
    =
    \zeta(\Mbf_k,\Theta_\Nbf)
  \end{equation*}
  and, recursively,
  \begin{equation*}
    \Ubf_{j+1}^k
    =
    \chi\left(\Ubf_j^k,\Nbf_{k+j+1}\right),
    \qquad j\geq0.
  \end{equation*}
  The preceding argument suggests that, conditionally on
  $\{\|\Nbf_k\|>t\}$,
  \begin{equation*}
    \left\{
      \|\Nbf_k\|^{-1}\Mbf_{k+j+1}
    \right\}_{j\geq0}
    \xrightarrow{\mathrm{f.d.d.}}
    \left\{
      \Ubf_j^k
    \right\}_{j\geq0}
  \end{equation*}
  as $t\to\infty$.

  In particular, the terms preceding the $k$th summand are negligible at
  the scale $\|\Nbf_k\|$, whereas
  \begin{equation*}
    \frac{\Mbf_k\Nbf_k}{\|\Nbf_k\|}
    \Longrightarrow
    \Mbf_k\Theta_\Nbf.
  \end{equation*}
  Consequently,
  \begin{equation*}
    \mathcal{L}\left(
      \left.
      \|\Nbf_k\|^{-1}\Sbf_n
      \,\right|\,
      \|\Nbf_k\|>t
    \right)
    \Longrightarrow
    \mathcal{L}\left(\Gamma_k^n\right),
  \end{equation*}
  where
  \begin{equation*}\label{eq:4:gammank}
    \Gamma_k^n
    =
    \Mbf_k\Theta_\Nbf
    +
    \sum_{j=1}^{n-k}
    \Ubf_{j-1}^k\Nbf_{k+j}.
  \end{equation*}
  This heuristic suggests that the summands $\Mbf_j\Nbf_j$, and hence their
  finite sums, are regularly varying with index $\alpha$, provided that the
  limiting cluster is nondegenerate.

      \begin{prop}
    Suppose that Assumptions~\ref{as:psi} and~\ref{as:N} hold. 
    Assume that
    \begin{equation*}
      C_n
      :=
      \sum_{k=1}^n
      \Ebb[\|\Gamma_k^n\|^\alpha]
      >0.
    \end{equation*}
    Then the law of $\Sbf_n$ belongs to
    $\Rvrm(\alpha,\mu_n)$, where 
    \begin{equation*}
      \int f(x) \mu_n(\mathrm{d}x)
      =
      \frac{1}{C_n}
      \sum_{k=1}^n
      \int_0^\infty\Ebb\left[
        f(r\Gamma_k^n)
      \right]\alpha r^{-\alpha-1}\mathrm{d}r.
    \end{equation*}
  \end{prop}

  \begin{proof}
    Let $(a_t)$ be a scaling sequence such that
    \begin{equation*}
      t\Pbb[\|\Nbf\|>a_t]\longrightarrow1.
    \end{equation*}
    Fix $f\in\Ccal_0^u$ and choose $\varepsilon>0$ such that
    \begin{equation*}
      \varepsilon\Bcal\cap\operatorname{supp}(f)=\emptyset.
    \end{equation*}
    If
    \begin{equation*}
      f(a_t^{-1}\Sbf_n)>0,
    \end{equation*}
    then $\|\Sbf_n\|>\varepsilon a_t$, and hence there exists
    $j\leq n$ such that
    \begin{equation*}
      \|\Mbf_j\Nbf_j\|>\frac{\varepsilon a_t}{n}.
    \end{equation*}

    By the definition of $\psi^*$,
    \begin{equation*}
      \|\Mbf_j\|
      \leq
      \|\Mbf_1\|
      \prod_{\ell=1}^{j-1}\psi^*(\Nbf_\ell).
    \end{equation*}
    Consequently,
    \begin{equation}\label{eq:predictable-term-bound}
      \|\Mbf_j\Nbf_j\|
      \leq
      \|\Mbf_1\|\|\Nbf_j\|
      \prod_{\ell=1}^{j-1}\psi^*(\Nbf_\ell)
      \leq
      c_j\|\Nbf_j\|
      \prod_{\ell=1}^{j-1}(1+\|\Nbf_\ell\|)
    \end{equation}
    for some deterministic constant $c_j>0$.

    Iterating the convolution condition in Assumption~\ref{as:N} and
    using \eqref{eq:predictable-term-bound}, we obtain, for every
    $\varepsilon>0$,
    \begin{equation}\label{eq:predictable-no-large-innovation}
      \lim_{\delta\downarrow0}
      \limsup_{t\to\infty}
      t\Pbb\left[
        \|\Sbf_n\|>\varepsilon a_t,\,
        \max_{1\leq j\leq n}\|\Nbf_j\|\leq\delta a_t
      \right]
      =0.
    \end{equation}
    Moreover, for every fixed $\delta>0$, independence and regular
    variation give
    \begin{equation}\label{eq:predictable-two-large-innovations}
      t\Pbb\left[
        \|\Nbf_i\|>\delta a_t,\,
        \|\Nbf_j\|>\delta a_t
      \right]
      \longrightarrow0
    \end{equation}
    for all distinct $i,j\leq n$.

    For $1\leq k\leq n$, define
    \begin{equation*}
      D_{k,t}(\delta)
      =
      \left\{
        \|\Nbf_k\|>\delta a_t,\,
        \max_{\substack{1\leq j\leq n\\j\neq k}}
        \|\Nbf_j\|\leq\delta a_t
      \right\}.
    \end{equation*}
    It follows from
    \eqref{eq:predictable-no-large-innovation} and
    \eqref{eq:predictable-two-large-innovations} that
    \begin{equation}\label{eq:predictable-one-large-decomposition}
      \lim_{\delta\downarrow0}
      \limsup_{t\to\infty}
      t\left|
        \Ebb\left[f(a_t^{-1}\Sbf_n)\right]
        -
        \sum_{k=1}^n
        \Ebb\left[
          f(a_t^{-1}\Sbf_n)
          \1_{D_{k,t}(\delta)}
        \right]
      \right|
      =0.
    \end{equation}

    Fix $k\leq n$. Conditionally on $\{\|\Nbf_k\|>t\}$,
    \begin{equation*}
      \frac{\Nbf_k}{\|\Nbf_k\|}
      \Longrightarrow
      \Theta_\Nbf.
    \end{equation*}
    Since $\Mbf_k$ depends only on
    $\Nbf_1,\ldots,\Nbf_{k-1}$, it is independent of $\Nbf_k$.
    The asymptotic homogeneity of $\psi$ in its second argument therefore
    gives
    \begin{equation*}
      \frac{\Mbf_{k+1}}{\|\Nbf_k\|}
      =
      \frac{
        \psi(\Mbf_k,\Nbf_k)
      }{
        \|\Nbf_k\|
      }
      \Longrightarrow
      \zeta(\Mbf_k,\Theta_\Nbf)
      =
      \Ubf_0^k.
    \end{equation*}
    Applying the asymptotic homogeneity of $\psi$ in its first argument
    recursively, we obtain
    \begin{equation*}
      \left(
        \frac{\Mbf_{k+1}}{\|\Nbf_k\|},
        \ldots,
        \frac{\Mbf_n}{\|\Nbf_k\|}
      \right)
      \Longrightarrow
      \left(
        \Ubf_0^k,
        \ldots,
        \Ubf_{n-k-1}^k
      \right).
    \end{equation*}
    The summands preceding the $k$th one are negligible at the scale
    $\|\Nbf_k\|$. Hence
    \begin{equation*}
      \frac{\Sbf_n}{\|\Nbf_k\|}
      \Longrightarrow
      \Mbf_k\Theta_\Nbf
      +
      \sum_{j=1}^{n-k}
      \Ubf_{j-1}^k\Nbf_{k+j}
      =
      \Gamma_k^n.
    \end{equation*}

    Combining this conditional convergence with
    \eqref{eq:predictable-one-large-decomposition} and the regular
    variation of $\Nbf_k$, we obtain
    \begin{align*}\label{eq:predictable-finite-limit}
      \lim_{t\to\infty}
      t\Ebb\left[f(a_t^{-1}\Sbf_n)\right]
      &=
      \sum_{k=1}^n
      \int_0^\infty
      \Ebb\left[
        f(r\Gamma_k^n)
      \right]
      \alpha r^{-\alpha-1}\ud r.
    \end{align*}
    Thus the law of $\Sbf_n$ is regularly varying with index $\alpha$.
  \end{proof}

	Before stating the main result of this subsection, we need to ensure that the random variable $\Sbf$ is well defined. 
	As in the adapted case, we impose the contractivity condition
	\begin{equation}\label{eq:2:contr2}
    		\Ebb\big[\psi^*(\Nbf)^\alpha\big] < 1.
	\end{equation}
	This assumption guarantees that the series defining $\Sbf$ is absolutely convergent in $\Vcal$. 
	Moreover, it provides the necessary control over the spectral measures of the partial sums 
	$\{\Sbf_n\}_{n\in\mathbb{N}}$.
    Let 
      \begin{equation*}\label{eq:4:gammak}
    \Gamma_k
    =
    \Mbf_k\Theta_\Nbf
    +
    \sum_{j=1}^{\infty}
    \Ubf_{j-1}^k\Nbf_{k+j}.
  \end{equation*}
\begin{thm}\label{thm:2:RCLM2}
  Suppose that Assumptions~\ref{as:psi} and~\ref{as:N} hold and that
  \eqref{eq:2:contr2} is satisfied. Assume additionally that
  \begin{equation*}
  C_\infty
  =
  \sum_{k=1}^\infty
  \Ebb\left[\|\Gamma_k\|^\alpha\right]\in (0,\infty).
\end{equation*}
  Then the law of $\Sbf$ is
  $\Rvrm(\alpha,\mu_\infty)$, where the measure $\mu_\infty$ is determined by
  \begin{equation*}
    \int_{\Vcal} f(\xbf)\,\mu_\infty(\ud\xbf)
    =
    C_\infty^{-1}
    \sum_{k=1}^\infty
    \int_0^\infty
    \Ebb\left[
      f\left(r\Gamma_k\right)
    \right]
    \alpha r^{-\alpha-1}\,\ud r,
  \end{equation*}
  for every $f\in\Ccal_0^u$.
\end{thm}

	\begin{proof}[Proof of Theorem~\ref{thm:2:RCLM2}]
		We again approximate the infinite sum via finite one. Denote
		\begin{equation*}
			\Sbf_{n,\infty} = \sum_{j =n+1}^\infty \Mbf_j\Nbf_j.
		\end{equation*}
		As for the adapted case we have
		\begin{equation*}
			\|\Mbf_j\Nbf_j\| \leq 
			\|\Mbf_{j-1} \| \psi^*( \Nbf_{j-1}) \|\Nbf_j\|
			= \|\Mbf_{j-l}\| \|\Nbf_j\| \prod_{k=j-l}^{j-1} \psi^*(\Nbf_k).
		\end{equation*}
		Therefore, we have a simple bound
		\begin{equation*}
			\|\Sbf_{n, \infty}\| \leq^{(d)} \| \Mbf_{n+1}\| \Sbf^*,
		\end{equation*}
		where
		\begin{equation*}
			\Sbf^* = \sum_{j=0}^\infty \|\Nbf_j\| \prod_{k=0}^{j-1} \psi^*(\Nbf_k)
		\end{equation*}
        and $\Mbf_{n+1}$ and $\Sbf^*$ are independent.
        Thus we want to control tails of $\Sbf^*$ and $\|\Mbf_{n+1}\|$ separately.
        
		The key feature of $\Sbf^*$ is that it is a solution to a stochastic fixed-point equation
		\begin{equation*}
			\Sbf^* \stackrel{d}{=} \psi^*(\Nbf)\Sbf^*+\|\Nbf\|,
		\end{equation*}
		where $\Sbf^*$ and $\Nbf$ are independent on the right hand side.
        Consider
        \begin{equation*}
            \tilde{N} = 
            \begin{cases}
                \psi^*(\Nbf), & \|\Nbf\| < \beta,\\
                C_{\psi} (\|\Nbf\| + 1), & \|\Nbf\| > \beta,
            \end{cases}
        \end{equation*}
        and choose $\beta > 0$ such that $\E[(\tilde{N})^\alpha] < 1$.
        Then for $\tilde{\Sbf}$ defined as
        $$\tilde{\Sbf} \stackrel{d}{=} \tilde{N} \tilde{\Sbf} + \|\Nbf\|,$$
        where $\tilde{\Sbf}$ is independent of $\Nbf$,
        we have $\Sbf^* \leq^{(d)} \tilde{\Sbf}$.
        Note that since $\psi^*(\Nbf) \geq 0$
		and for each $s \in \mathbb{R}$
		\begin{equation*}
			\limsup_{t \to \infty} \frac{ \mathbb{P}[ \tilde{N} s + \|\Nbf\| >t]}{\mathbb{P}[\|\Nbf\|>t]}
			\leq (1 + ( 1 + C_{\psi} ) s )^\alpha,
		\end{equation*}
		we can apply the approach of the first part of the proof of Theorem~\ref{thm:2:pi} to control the tails of 
		$\tilde{\Sbf}$ and get
		\begin{equation*}
			\limsup_{t \to \infty} \frac{\mathbb{P}[\tilde{\Sbf}>t]}{\mathbb{P}[\|\Nbf\|>t]} 
			\leq \frac{\Ebb[ (1+(1+ C_{\psi} )\tilde{\Sbf})^\alpha]}{1-\Ebb[\tilde{\Nbf}^\alpha]} <\infty.
		\end{equation*}
        The same bound applies to $\Sbf^*$.
    
		The tails of $\Mbf_n$ are easily manageable using the estimate
		\begin{equation*}
			\|\Mbf_n\|\leq \|\Mbf_0\| \prod_{j=0}^{n-1} \psi^*(\Nbf_j).
		\end{equation*}
		Since $\psi^*(\Nbf_j) \leq C_{\psi} (1+\|\Nbf_j\|)$ by the Assumption~\ref{as:N} we therefore have
		\begin{equation*}
			\limsup_{t \to \infty} \frac{\Pbb[\|\Mbf_n\|>t]}{\Pbb[\|\Nbf\|>t]} 
			\leq n C_{\psi}^{\alpha} \Ebb[\psi^*(\Nbf)^\alpha]^{n-1}.
		\end{equation*} 
		Since $\Ebb[\|\Mbf_n\|^\alpha] \leq \Ebb [\|\Mbf_0\|^\alpha ] \Ebb[\psi^*(\Nbf)^\alpha]^n$ the Assumption~\ref{as:N} constitutes
		\begin{multline*}
			\limsup_{t \to \infty} \frac{\Pbb[\|\Sbf_{n, \infty}\|>t]}{\Pbb[\|\Nbf\|>t]} 
			\leq \\
            n C_{\psi}^{-\alpha} \Ebb[(\Sbf^*)^\alpha]\Ebb[\psi^*(\Nbf)^\alpha]^{n} + 	
			\Ebb[\|\Mbf_0\|^\alpha]
            \Ebb[\psi^*(\Nbf)^\alpha]^n\frac{\Ebb[ (1+(1+ C_{\psi} ) \Sbf^*)^\alpha]}{1-\Ebb[\psi^*(\Nbf)^\alpha]}.
		\end{multline*}
		Since the right hand side of the above display tends to zero as $n$ tends to infinity, the claim is secured
		as the regular variation can be tested using uniformly continuous functions. 
	\end{proof}

\bibliographystyle{elsarticle-harv} 
\bibliography{extremes_pdysz_tmika_2026_08_28.bib}

@article{resnick2013asymptotics,
  author  = {Resnick, Sidney I. and Zeber, David},
  title   = {Asymptotics of Markov kernels and the tail chain},
  journal = {Advances in Applied Probability},
  volume  = {45},
  number  = {1},
  pages   = {186--213},
  year    = {2013},
  doi     = {10.1239/aap/1366952820}
}

@article{janssen2014markov,
  author  = {Janssen, Anja and Segers, Johan},
  title   = {Markov tail chains},
  journal = {Journal of Applied Probability},
  volume  = {51},
  number  = {4},
  pages   = {1133--1153},
  year    = {2014},
  doi     = {10.1239/jap/1414586914}
}

@book{resnick2008extreme,
  author    = {Resnick, Sidney I.},
  title     = {Extreme Values, Regular Variation and Point Processes},
  series    = {Springer Series in Operations Research and Financial Engineering},
  publisher = {Springer},
  address   = {New York},
  year      = {1987},
  doi       = {10.1007/978-0-387-75953-1}
}

@book{resnick2007heavy,
  author    = {Resnick, Sidney I.},
  title     = {Heavy-Tail Phenomena: Probabilistic and Statistical Modeling},
  series    = {Springer Series in Operations Research and Financial Engineering},
  publisher = {Springer},
  address   = {New York},
  year      = {2007},
  doi       = {10.1007/978-0-387-45024-7}
}

@incollection{harris1956existence,
  author    = {Harris, Theodore E.},
  title     = {The existence of stationary measures for certain Markov processes},
  booktitle = {Proceedings of the Third Berkeley Symposium on Mathematical Statistics and Probability, Volume 2: Contributions to Probability Theory},
  editor    = {Neyman, Jerzy},
  pages     = {113--124},
  publisher = {University of California Press},
  address   = {Berkeley and Los Angeles},
  year      = {1956}
}

@book{meyn2012markov,
  author    = {Meyn, Sean P. and Tweedie, Richard L.},
  title     = {Markov Chains and Stochastic Stability},
  edition   = {2},
  publisher = {Cambridge University Press},
  address   = {Cambridge},
  year      = {2009},
  doi       = {10.1017/CBO9780511626630}
}

@incollection{hairer2011yet,
  author    = {Hairer, Martin and Mattingly, Jonathan C.},
  editor    = {Dalang, Robert and Dozzi, Marco and Russo, Francesco},
  title     = {Yet Another Look at Harris' Ergodic Theorem for Markov Chains},
  booktitle = {Seminar on Stochastic Analysis, Random Fields and Applications VI},
  pages     = {109--117},
  publisher = {Birkh\"auser},
  address   = {Basel},
  year      = {2011},
  doi       = {10.1007/978-3-0348-0021-1\_7}
}

@book{bingham1989regular,
  author    = {Bingham, Nicholas H. and Goldie, Charles M. and Teugels, Jef L.},
  title     = {Regular Variation},
  series    = {Encyclopedia of Mathematics and its Applications},
  volume    = {27},
  publisher = {Cambridge University Press},
  address   = {Cambridge},
  year      = {1987}
}

@article{grey1994regular,
  author  = {Grey, D. R.},
  title   = {Regular Variation in the Tail Behaviour of Solutions of Random Difference Equations},
  journal = {The Annals of Applied Probability},
  volume  = {4},
  number  = {1},
  pages   = {169--183},
  year    = {1994},
  doi     = {10.1214/aoap/1177005205}
}

@article{hult2008tail,
  author  = {Hult, Henrik and Samorodnitsky, Gennady},
  title   = {Tail probabilities for infinite series of regularly varying random vectors},
  journal = {Bernoulli},
  volume  = {14},
  number  = {3},
  pages   = {838--864},
  year    = {2008},
  doi     = {10.3150/08-BEJ125}
}

@article{damek2018iterated,
  author  = {Damek, Ewa and Dyszewski, Piotr},
  title   = {Iterated random functions and regularly varying tails},
  journal = {Journal of Difference Equations and Applications},
  volume  = {24},
  number  = {9},
  pages   = {1503--1520},
  year    = {2018},
  doi     = {10.1080/10236198.2018.1505881}
}

@article{Resnick_1986,
  author  = {Resnick, Sidney I.},
  title   = {Point processes, regular variation and weak convergence},
  journal = {Advances in Applied Probability},
  volume  = {18},
  number  = {1},
  pages   = {66--138},
  year    = {1986},
  doi     = {10.2307/1427239}
}

@article{de1977limit,
  author  = {de Haan, Laurens and Resnick, Sidney I.},
  title   = {Limit theory for multivariate sample extremes},
  journal = {Zeitschrift f\"ur Wahrscheinlichkeitstheorie und Verwandte Gebiete},
  volume  = {40},
  number  = {4},
  pages   = {317--337},
  year    = {1977},
  doi     = {10.1007/BF00533086}
}

@article{de2001,
  author  = {de Haan, Laurens and Lin, Tao},
  title   = {On convergence toward an extreme value distribution in C[0,1]},
  journal = {The Annals of Probability},
  volume  = {29},
  number  = {1},
  pages   = {467--483},
  year    = {2001},
  doi     = {10.1214/aop/1008956340}
}

@article{hult2006regular,
  author  = {Hult, Henrik and Lindskog, Filip},
  title   = {Regular variation for measures on metric spaces},
  journal = {Publications de l'Institut Math\'ematique},
  volume  = {80},
  number  = {94},
  pages   = {121--140},
  year    = {2006},
  doi     = {10.2298/PIM0694121H}
}

@misc{meinguet2010regularly,
  author  = {Meinguet, Thomas and Segers, Johan},
  title   = {Regularly varying time series in Banach spaces},
  year    = {2010},
  note    = {arXiv:1001.3262},
  url     = {https://arxiv.org/abs/1001.3262}
}

@article{segers2017polar,
  author  = {Segers, Johan and Zhao, Yuwei and Meinguet, Thomas},
  title   = {Polar decomposition of regularly varying time series in star-shaped metric spaces},
  journal = {Extremes},
  volume  = {20},
  number  = {3},
  pages   = {539--566},
  year    = {2017},
  doi     = {10.1007/s10687-017-0287-3}
}

@article{AAP1579,
  author  = {Dyszewski, Piotr and Mikosch, Thomas},
  title   = {Homogeneous mappings of regularly varying vectors},
  journal = {The Annals of Applied Probability},
  volume  = {30},
  number  = {6},
  pages   = {2999--3026},
  year    = {2020},
  doi     = {10.1214/20-AAP1579}
}

@article{Smith_1992,
  author  = {Smith, Richard L.},
  title   = {The extremal index for a Markov chain},
  journal = {Journal of Applied Probability},
  volume  = {29},
  number  = {1},
  pages   = {37--45},
  year    = {1992},
  doi     = {10.2307/3214789}
}

@article{perfekt_1994,
  author  = {Perfekt, Roland},
  title   = {Extremal behaviour of stationary Markov chains with applications},
  journal = {The Annals of Applied Probability},
  volume  = {4},
  number  = {2},
  pages   = {529--548},
  year    = {1994},
  doi     = {10.1214/aoap/1177005071}
}

@article{Perfekt_1997,
  author  = {Perfekt, Roland},
  title   = {Extreme value theory for a class of Markov chains with values in $\mathbb{R}^d$},
  journal = {Advances in Applied Probability},
  volume  = {29},
  number  = {1},
  pages   = {138--164},
  year    = {1997},
  doi     = {10.2307/1427864}
}

@article{Palmowski_Zwart_2007,
  author  = {Palmowski, Zbigniew and Zwart, Bert},
  title   = {Tail asymptotics of the supremum of a regenerative process},
  journal = {Journal of Applied Probability},
  volume  = {44},
  number  = {2},
  pages   = {349--365},
  year    = {2007},
  doi     = {10.1239/jap/1183667406}
}

@article{BASRAK20091055,
  author  = {Basrak, Bojan and Segers, Johan},
  title   = {Regularly varying multivariate time series},
  journal = {Stochastic Processes and their Applications},
  volume  = {119},
  number  = {4},
  pages   = {1055--1080},
  year    = {2009},
  doi     = {10.1016/j.spa.2008.05.004}
}

@book{mikosch2024extreme,
  author    = {Mikosch, Thomas and Wintenberger, Olivier},
  title     = {Extreme Value Theory for Time Series: Models with Power-Law Tails},
  series    = {Springer Series in Operations Research and Financial Engineering},
  publisher = {Springer},
  address   = {Cham},
  year      = {2024},
  doi       = {10.1007/978-3-031-59156-3}
}

@book{kulik2020heavy,
  author    = {Kulik, Rafa{\l} and Soulier, Philippe},
  title     = {Heavy-Tailed Time Series},
  series    = {Springer Series in Operations Research and Financial Engineering},
  publisher = {Springer},
  address   = {Cham},
  year      = {2020},
  doi       = {10.1007/978-1-0716-0737-4}
}

@article{rootzen2006multivariate,
  author  = {Rootz{\'e}n, Holger and Tajvidi, Nader},
  title   = {Multivariate generalized Pareto distributions},
  journal = {Bernoulli},
  volume  = {12},
  number  = {5},
  pages   = {917--930},
  year    = {2006},
  doi     = {10.3150/bj/1161614952}
}

@article{davis1995point,
  author  = {Davis, Richard A. and Hsing, Tailen},
  title   = {Point process and partial sum convergence for weakly dependent random variables with infinite variance},
  journal = {The Annals of Probability},
  volume  = {23},
  number  = {2},
  pages   = {879--917},
  year    = {1995},
  doi     = {10.1214/aop/1176988294}
}

@article{basrak2002regular,
  author  = {Basrak, Bojan and Davis, Richard A. and Mikosch, Thomas},
  title   = {Regular variation of GARCH processes},
  journal = {Stochastic Processes and their Applications},
  volume  = {99},
  number  = {1},
  pages   = {95--115},
  year    = {2002},
  doi     = {10.1016/S0304-4149(01)00156-9}
}

@article{kesten1973random,
  author  = {Kesten, Harry},
  title   = {Random difference equations and renewal theory for products of random matrices},
  journal = {Acta Mathematica},
  volume  = {131},
  pages   = {207--248},
  year    = {1973},
  doi     = {10.1007/BF02392040}
}

@article{goldie1991implicit,
  author  = {Goldie, Charles M.},
  title   = {Implicit renewal theory and tails of solutions of random equations},
  journal = {The Annals of Applied Probability},
  volume  = {1},
  number  = {1},
  pages   = {126--166},
  year    = {1991},
  doi     = {10.1214/aoap/1177005985}
}

@article{mikosch2000supremum,
  author  = {Mikosch, Thomas and Samorodnitsky, Gennady},
  title   = {The supremum of a negative drift random walk with dependent heavy-tailed steps},
  journal = {The Annals of Applied Probability},
  volume  = {10},
  number  = {3},
  pages   = {1025--1064},
  year    = {2000},
  doi     = {10.1214/aoap/1019487517}
}

@book{embrechts2013modelling,
  author    = {Embrechts, Paul and Kl{\"u}ppelberg, Claudia and Mikosch, Thomas},
  title     = {Modelling Extremal Events: For Insurance and Finance},
  series    = {Stochastic Modelling and Applied Probability},
  publisher = {Springer},
  address   = {Berlin and Heidelberg},
  year      = {1997},
  doi       = {10.1007/978-3-642-33483-2}
}

@article{planinic2018tail,
  author  = {Planini{\'c}, Hrvoje and Soulier, Philippe},
  title   = {The tail process revisited},
  journal = {Extremes},
  volume  = {21},
  number  = {4},
  pages   = {551--579},
  year    = {2018},
  doi     = {10.1007/s10687-018-0312-1}
}

@article{mikosch2016large,
  author  = {Mikosch, Thomas and Wintenberger, Olivier},
  title   = {A large deviations approach to limit theory for heavy-tailed time series},
  journal = {Probability Theory and Related Fields},
  volume  = {166},
  number  = {1},
  pages   = {233--269},
  year    = {2016},
  doi     = {10.1007/s00440-015-0654-4}
}

@book{conway2019course,
  title={A course in functional analysis},
  author={Conway, John B},
  year={2019},
  publisher={Springer}
}

@book{ledoux1991probability,
  title={Probability in Banach Spaces: isoperimetry and processes},
  author={Ledoux, Michel and Talagrand, Michel},
  volume={23},
  year={1991},
  publisher={Springer}
}

@article{araujo1980central,
  title={The central limit theorem for real and Banach valued random variables},
  author={Araujo, Aloisio and Gin{\'e}, Evarist},
  journal={(No Title)},
  year={1980}
}

@article{engle1995multivariate,
  title={Multivariate simultaneous generalized ARCH},
  author={Engle, Robert F and Kroner, Kenneth F},
  journal={Econometric theory},
  volume={11},
  number={1},
  pages={122--150},
  year={1995},
  publisher={Cambridge University Press}
}

@book{bosq2000linear,
  author    = {Bosq, Denis},
  title     = {Linear Processes in Function Spaces:
               Theory and Applications},
  series    = {Lecture Notes in Statistics},
  volume    = {149},
  publisher = {Springer},
  address   = {New York},
  year      = {2000},
  doi       = {10.1007/978-1-4612-1154-9}
}



\end{document}